\documentclass{amsart}
\usepackage{amsmath}
\usepackage{amsthm}
\usepackage{amssymb}

 \newtheorem{thm}{Theorem}[subsection]
 \newtheorem{cor}[thm]{Corollary}
 \newtheorem{lem}[thm]{Lemma}
 \newtheorem{prop}[thm]{Proposition}
 \theoremstyle{definition}
 
 \theoremstyle{remark}
 \newtheorem{rem}[thm]{Remark}
 \numberwithin{equation}{subsection}
\newcommand{\pn}{\noindent}

\newcommand{\ZZ}{\mathbb{Z}}

\newcommand{\LL}{\mathbb{L}}

\newcommand{\GG}{\mathbb{G}}
\newcommand{\PP}{\mathbb{P}}

\newcommand{\li}{\mathrm{Lie}\,}

\newcommand{\Hom}{\mathrm{Hom}}
\newcommand{\Ext}{\mathrm{Ext}}
\newcommand{\Biext}{\mathrm{Biext}}
\newcommand{\uHom}{\underline{\mathrm{Hom}}}
\newcommand{\uExt}{\underline{\mathrm{Ext}}}
\newcommand{\uMor}{\underline{\mathrm{Mor}}}
\newcommand{\bBiext}{\mathrm{\mathbf{Biext}}}

\newcommand{\bExt}{\mathrm{\mathbf{Ext}}}
\newcommand{\W}{\mathrm{W}}
\newcommand{\F}{\mathrm{F}}

\newcommand{\h}{\mathrm{H}}
\newcommand{\R}{\mathrm{R}}

\newcommand{\rk}{\mathrm{rk}}

\newcommand{\spec}{{\mathrm{Spec}}\,}

\newcommand{\oks}{\overline {k(s)}}
\newcommand{\os}{\overline s}

\begin{document}

\title[extensions and biextensions]
{Extensions and biextensions of locally constant group schemes, tori and abelian schemes}

\author{Cristiana Bertolin}

\address{NWF-I Mathematik, Universit\"at Regensburg, D-93040 Regensburg}

\email{cristiana.bertolin@mathematik.uni-regensburg.de}

\subjclass{18A20;14A15}

\keywords{extensions, biextensions, locally constant group schemes, tori, abelian schemes}

\date{}
\dedicatory{}

\commby{Cristiana Bertolin}

%%% ----------------------------------------------------------------------

\begin{abstract}
Let $S$ be a scheme. We compute explicitly the group of homomorphisms, the $S$-sheaf of homomorphisms, the group of extensions, and the $S$-sheaf of extensions involving locally constant $S$-group schemes, abelian $S$-schemes, and $S$-tori. Using the obtained results, we study the categories of biextensions involving these geo\-metrical objets. In particular, we prove that if $G_i$ (for $i=1,2,3$) is an extension of an abelian
$S$-scheme $A_i$ by an $S$-torus $T_i$, the category of biextensions of $(G_1,G_2)$ by $G_3$ is equivalent to the category of biextensions of the underlying abelian $S$-schemes $(A_1,A_2)$ by the underlying $S$-torus $T_3$.

\end{abstract}

%%% ----------------------------------------------------------------------
\maketitle
%%% ----------------------------------------------------------------------

%\tableofcontents

\section*{Introduction}

The notion of biextension was introduced by D. Mumford~\cite{Mu1} in the context of formal groups in order to express the relations between the formal groups associated to an abelian scheme and the dual abelian scheme. Successively in the Expos\'es VII and VIII of~\cite{SGA7}, Grothendieck studies in a systematic way the notion of biextension in the more general setting of abelian sheaves over any topos. In particular he investigates the case of biextensions of commutative group schemes by the multiplicative group ${\GG}_m$.

The aim of this paper is to study biextensions involving locally constant group schemes, tori and abelian schemes defined over an arbitrary base scheme $S$.
We treat all possible cases - biextensions of abelian schemes by abelians schemes, biextensions of tori and abelian schemes by locally constant schemes, ... - and in order to have a self-contained work, we recall briefly the cases which already exist in the literature. Working in the topos $\mathrm{{\mathbf{T}_{fppf}}}$ associated to the site of locally of finite presentation $S$-schemes  endowed with the fppf topology, our main results are (respectively Theorems~\ref{biextPTA},~\ref{biextAAA},~\ref{mthm1}):

 \textbf{Theorem A}:\emph{
 Let $A$ be an abelian $S$-scheme, let $T$ be an $S$-torus and let $P$ be a divisible commutative $S$-group scheme locally of finite presentation over $S$, with connected fibres.
The category $ {\bBiext}(P,T;A)$ of biextensions of $(P,T)$ by $A$
and the category $ {\bBiext}(T,P;A)$ of biextensions of $(T,P)$ by $A$ are equivalent to the trivial category.}

 \textbf{Theorem B}:\emph{
Let $A_i$ (for $i=1,2,3$) be an abelian $S$-scheme.
The category\\
 $ {\bBiext}(A_1,A_2;A_3)$ of biextensions of $(A_1,A_2)$ by $A_3$ is equivalent to the trivial category.}

 \textbf{Theorem C}:\emph{
Let $S$ be a scheme.
Let $G_i$ (for $i=1,2,3$) be a commutative extension of an abelian
$S$-scheme $A_i$ by an $S$-torus $T_i$. The category of biextensions of $(G_1,G_2)$ by $G_3$ is equivalent
to the category of biextensions of the underlying abelian $S$-schemes $(A_1,A_2)$ by the underlying $S$-torus $T_3$.}\\
The reason of the choice of the topos $\mathrm{{\mathbf{T}_{fppf}}}$ is that
for the fppf topology the tori and the abelian schemes are divisible groups.

The theory of motives leads us to the investigation of biextensions involving locally constant group schemes, tori and abelian schemes defined over an arbitrary base scheme $S$. In fact in~\cite{B} we introduce the notion of biextension of 1-motives by 1-motives and we define bilinear morphisms between 1-motives as isomorphism classes of such biextensions. We then check that our definition is compatible with the realizations of 1-motives, i.e. that the group of isomorphism classes of biextensions of 1-motives is a group of bilinear morphisms in an appropriate category of mixed realizations.
In this context \textbf{Theorem A}, \textbf{B}
 and \textbf{C} mean that biextensions satisfy the main property of morphisms of motives, i.e. they respect the weight filtration ${\W}_*$ on motives.
  For example, if $A_i$ (for $i=1,2,3$) is an abelian $S$-scheme,
\textbf{Theorem B} says that there are no nonzero bilinear morphisms from $A_1 \times A_2$ to $A_3$. But this is exactly what is predicted by Grothendieck's philosophy of motives: since morphisms of motives have to respect the weight filtration ${\W}_*$, it is not possible to have a morphism from the motive $A_1 \otimes A_2$ of weight -2 to the motive $A_3$ of weight -1.

Before we investigate
biextensions, we have to understand extensions and homomorphisms: working over an arbitrary base scheme $S$, we start  computing explicitly the group of homomorphisms, the $S$-sheaf of homomorphisms, the group of extensions and the $S$-sheaf of extensions involving locally constant $S$-group schemes, abelian $S$-schemes and $S$-tori. There are a lot of results in the literature about homomorphisms between such geometric objects but very few about their extensions.
This relies on the fact that with extensions of group schemes we go outside the category of schemes: extensions of $S$-group schemes are a priori only algebraic spaces over $S$ and this makes things immediately more complicated since in the literature there are only few results about algebraic spaces that are relevant to the study of extensions considered in this paper.  This will be the most difficult point in the study of extensions of an $S$-torus $T$ by an abelian $S$-scheme $A$, which begins with the proof that
 for such extensions it is equivalent to be a scheme or to be of finite order locally over $S$ (Proposition~\ref{extTA-equiv}). Using the fact that
these extensions are representable by schemes (Theorem~\ref{extTA}), we can
then conclude that the extensions of $S$-tori by abelian $S$-schemes are of finite order locally over $S$, i.e. the $S$-sheaf ${\uExt}^1(T,A)$ is a torsion sheaf (Corollary~\ref{extTA-cor1}). Moreover we show that this $S$-sheaf is in fact an $S$-group scheme which is separated and \'etale over $S$ (Corollary~\ref{extTA-cor3}). With extensions of abelian $S$-schemes by abelian $S$-schemes we don't have problems of representability by schemes since it is a classical result that abelian algebraic spaces are abelian schemes. The main difficulty with these extensions is that they are not of finite order locally over $S$ (it is true only if we assume the base scheme $S$ to be integral and geometrically unibranched, see Proposition~\ref{extAAtorsion}) and we cannot say much about them.

 Because of the homological interpretation of biextensions furnished by Gro\-then\-dieck in~\cite{SGA7} Expos\'e VII 3.6.5 and (3.7.4), our results on biextensions are essentially consequences of the results on extensions obtained in the first part of this paper. Through a homological ``d\'evissage'', the fact that the $S$-sheaf ${\uExt}^1(T,A)$ is a separated and \'etale $S$-group scheme implies \textbf{Theorem A}.
The proof of \textbf{Theorem B} is done again trough a homological ``d\'evissage''
but it is not based on results of the first part of this paper since, as we have already said, extensions of abelian schemes by abelian schemes are not of finite order locally over $S$ if $S$ fails to be integral and geometrically unibranched!
 Concerning \textbf{Theorem C}, in~\cite{SGA7} Expos\'e VIII (3.6.1) Gro\-then\-dieck proves that in the topos
$\mathrm{{\mathbf{T}_{fppf}}}$, the category of biextensions of $(G_1,G_2)$ by ${\GG}_m$ is equivalent to the category of biextensions of the underlying abelian $S$-schemes $(A_1,A_2)$ by ${\GG}_m$. Therefore the proof of \textbf{Theorem C}
reduces to verifying that the category of biextensions of $(G_1,G_2)$ by  $G_3$ is equivalent to the category of biextensions of $(G_1,G_2)$ by the underlying torus $T_3.$
Because of the homological interpretation of biextensions, the proof
of this last equivalence of categories is an easy consequence of  \textbf{Theorem A} and \textbf{B}.

This article is a shortened version of the first two chapters of author's Habilitationsschrift~\cite{B}.

Je tiens \`a remercier L. Illusie, L. Moret-Bailly et M. Raynaud de leur aide pr\'ecieuse lors de l'\'elaboration des parties techniques de ce papier.
I am grateful also to the anonymous referee for the useful remarks.

%------------------------------------------------------------------------

\section*{Notation}

In this paper $S$ is an arbitrary scheme.

Let $s$ a point of $S$, with residue field $k(s)$. We denote by
 $\os$ a geometric point over $s$ and by $\oks$ its residue field.

If $P$, $Q$ are $S$-group schemes, we write $P_Q$ for the fibred product $P \times_S Q$ of $P$ and $Q$ over $S$, viewed as scheme over $Q$. In particular, if $s$ is a point of $S$, we denote by $P_s= P \times_S {\spec} (k(s))$ the fibre of $P$ over $s$.

An \textbf{abelian $S$-scheme} is an $S$-group scheme which is smooth, proper over $S$ and with connected fibres. An \textbf{$S$-torus} is an $S$-group scheme which is locally isomorphic for the fpqc topology (equivalently for the \'etale topology) to a $S$-group scheme of the kind ${\GG}_m^r$ with $r$ an integer bigger or equal to 0.
The \textbf{character group} ${\uHom}(T,{\GG}_m)$ and the
\textbf{cocharacter group} $ {\uHom}({\GG}_m,T)$
of a $S$-torus $T$
are $S$-group schemes which are \textbf{locally}
for the \'etale topology \textbf{constant group schemes defined by
 finitely generated free $\ZZ$-modules}. These $S$-group schemes are just locally of finite presentation over $S$ (not of finite presentation), and so in some later considerations it will be necessary to allow $S$-group schemes that are merely locally of finite presentation over $S$.
Sometimes it will be more convenient to denote by $Y(1)$ a torus with cocharacter group $Y={\uHom}({\GG}_m,Y(1))$ and character group $Y^{\vee}={\uHom}(Y(1),{\GG}_m)$.

 An $S$-group scheme is \textbf{divisible} if for each integer $n, n \not=0,$ the ``multiplication by $n$'' on it is an epimorphism for the fppf topology.

An \textbf{isogeny} $f:P \rightarrow Q$ between $S$-group schemes
is a morphism of $S$-group schemes which is finite, faithfully flat and of finite presentation.

An \textbf{algebraic space over $S$} is a functor $X:(S-Schemes)^\circ \rightarrow (Sets)$ satisfying the following conditions:
\begin{enumerate}
    \item $X$ is a sheaf for the \'etale topology,
    \item $X$ is locally representable: there exists an $S$-scheme $U$, and a map $U \rightarrow X$ which is representable by \'etale surjective maps,
    \item $X$ is quasi-separated over $S$.
\end{enumerate}

% ------------------------------------------------------------------------

\section{Homomorphisms and extensions}

Let $S$ be a scheme.

We are interested only in commutative extensions of commutative $S$-group sche\-mes. We consider such extensions in the category of abelian $S$-sheaves for the fppf site over $S$.

\begin{lem}\label{algspace}
Any commutative extension of commutative $S$-group schemes is an algebraic space over $S$ which is a group object.
\end{lem}

\begin{proof} Since by~\cite{SGA7} Expos\'e VII 1.2 a commutative extension of commutative $S$-group schemes is a torsor endowed with a group law, this Lemma is a consequence of~\cite{LM-B} Corollary 10.4.2.
\end{proof}

\subsection{The case in which locally constant group schemes are involved}

We start by studying the homomorphisms from a group scheme with connected fibres to an separated and unramified group scheme:

\begin{lem}\label{conetal}
Let $X$ be a commutative $S$-group scheme separated and unramified over $S$, and let $P$ be a commutative $S$-group scheme locally of finite presentation over $S$, with connected fibres. Then
\[{\Hom}(P,X)=0 .\]
\end{lem}

\begin{proof} Let $f: P \rightarrow X$ be an $S$-homomorphism and consider the $S$-morphism
$(f, \epsilon \circ p)_S: P \rightarrow X \times_{S} X,$
where $ \epsilon: S \rightarrow X$ is the unit section of $X$ and $p:P \rightarrow S$ is the structural morphism of $P$. According to~\cite{EGAIV} 4 Proposition 17.4.6, the inverse image
$(f, \epsilon \circ p)_S^{-1}(\Delta_{X/S})$ of the diagonal
 $ \Delta_{X/S}$ of $X$
 is an open and closed subscheme of $P$ whose restriction
 over each point of $S$ is not empty. Since
 the fibres of $P$ are connected, this inverse image
$(f, \epsilon \circ p)_S^{-1}(\Delta_{X/S})$
 is equal to $P$ and so $f$ is trivial.
\end{proof}

Concerning extensions of group schemes by locally constant group schemes, by~\cite{SGA3} Expos\'e X 5.5 we have that

\begin{lem}
Let $X$ be an $S$-group scheme which is locally
for the \'etale topology a constant group scheme defined by a
 finitely generated free $\ZZ$-module. Any $X$-torsor is a scheme.
 In particular, any extension of an $S$-group scheme by $X$ is a scheme.
\end{lem}

In \cite{SGA7} Expos\'e VIII Proposition 3.4 Grothendieck shows that there are no non-trivial extensions of smooth $S$-group schemes with connected fibres by $S$-group schemes which are locally
for the \'etale topology constant group schemes defined by
 finitely generated free $\ZZ$-modules. Here we prove more in general that

\begin{prop}\label{homextGX}
Let $X$ be commutative $S$-group scheme separated, unramified over $S$ and with constant geometric fibres defined by
 finitely generated free $\ZZ$-modules, and let $P$ be a smooth commutative $S$-group scheme with connected fibres. Then the category ${\bExt}(P,X)$ of extensions of $P$ by
$X$ is the trivial category. More precisely,
\[{\Hom}(P,X)=0 \qquad \textrm{and}  \qquad {\Ext}^1(P,X)=0.\]
In particular, both $S$-sheaves ${\uHom}(P,X)$ and ${\uExt}^1(P,X)$ are trivial.
\end{prop}

\begin{proof} Let $E$ be an extension of $P$ by
$X$. We first show that this extension $E$ is an $S$-group scheme. Assume $S$ to be affine and Noetherian.
According to Lemma~\ref{algspace} the extension $E$ is an algebraic space over $S$ which is a group object.
By hypothesis the $X$-torsor
$E$ is separated, unramified (hence locally of finite presentation, locally quasi-finite) over $P$ and so according to~\cite{A68} Theorem 3.3 (or~\cite{K} II Corollary 6.16) $E$ is a scheme, which is separated and unramified over $P$.\\
By Lemma~\ref{conetal} it is enough to prove that the extension $E$ is trivial locally for the Zariski topology, i.e. we can suppose $S$ to be the spectrum of a local Artin ring. Then by~\cite{SGA3} Expos\'e $\mathrm{VI_A}$ 2.3 and Proposition 2.4, the connected component of the identity of the extension $E$ exists. We denote it by $E^\circ $.
The scheme $P$ is of finite presentation over the affine Noetherian scheme $S$. Therefore $P$ is also Noetherian and in particular it is a finite disjoint union of its connected components. Restricting over one of these components we can suppose $P$ connected.
Since the structural morphism
$p: E \rightarrow P$ of the torsor $E$ over $P$ is separated and unramified, according to~\cite{EGAIV} 4 Proposition (17.4.9), to the connected component
$E^\circ $ corresponds a section $s:P \rightarrow E$ of $p$
which is an isomorphism from $P$ to $ E^\circ$ and therefore the extension $E$ is trivial.
\end{proof}

According to~\cite{SGA7} Expos\'e VII 1.4 for any commutative $S$-group scheme $P$ the group ${\Ext}^1({\ZZ},P)$ is isomorphic to the group of isomorphism classes of $P$-torsors with respect to the fppf topology. Remark that if $P$ is smooth over $S$, $P$-torsors for the fppf topology are the same as $P$-torsors for the \'etale topology. We have the following lemma

\begin{lem}\label{uextXP}
Let $X$ be a $S$-group scheme which is locally
for the \'etale topology a constant group scheme defined by a
 finitely generated free $\ZZ$-module and let $P$ be any
 commutative $S$-group scheme.
The $S$-sheaf ${\uExt}^1(X,P)$ is trivial.
\end{lem}

\subsection{The case in which tori are involved}
We start investigating homomorphisms between an $S$-torus $T$ and an abelian $S$-scheme $A$. In~\cite{SGA7} Expos\'e VII 1.3.8 Gro\-then\-dieck proves that
\[{\uHom}(A,T)=0. \]

\begin{lem}\label{homTA}
Let $A$ be an abelian $S$-scheme and let $T$ be
an $S$-torus. Then
\[{\Hom}(T,A)=0. \]
In particular, the $S$-sheaf ${\uHom}(T,A)$ is trivial.
\end{lem}

\begin{proof} Let $f:T \rightarrow A$ be a morphism from $T$ to $A$.
 Consider the restriction $f_{i}:T[q^i] \rightarrow A[q^i]$ of $f$ to the points of order $q^i$ with $q$ an integer bigger than 1 and invertible over $S$, and $i>0.$ For each $i$, $f_i$ is a morphism of finite and \'etale $S$-schemes which factors set-theoretically through the unit section of $A[q^i]$,
 and so it is trivial. Since
the family $(T[q^i])_i$ is schematically dense (see~\cite{SGA3} Expos\'e IX \S 4), we can conclude.
\end{proof}

We now investigate the extensions involving abelians schemes and tori.
 Since $S$-tori are affine over $S$, according to~\cite{SGA1} Expos\'e VIII 2.1 we have that

\begin{lem}
Any torsor under an $S$-torus is a scheme. In particular, any extension of an $S$-group scheme by an $S$-torus is a scheme.
\end{lem}

It is a classical result that over a separably closed field extensions of tori by tori are trivial. This is no longer true over an arbitrary scheme $S$. Nevertheless,
in~\cite{SGA7} Expos\'e VIII Proposition 3.3.1 Grothendieck proves that
if $T_i$ (for $i=1,2$) is
a torus over an arbitrary scheme $S$, then
\[{\uExt}^1(T_1,T_2)=0.\]

Over an algebraically closed field $k$ the extensions of an abelian $k$-variety $A$ by the $k$-torus ${\GG}_m$ are far from trivial: they are parametrized by the $k$-rational points of the dual abelian variety of $A$. More in general,
if $A$ is an abelian scheme over an arbitrary scheme $S$, we have that
\[{\uExt}^1(A,{\GG}_m) = A^*.\]
where $A^*$ is the dual abelian scheme of $A$ (\cite{O} Chapter I \S 5).

In the literature there are only few results about extensions of group schemes by abelian schemes (see~\cite{S} 7.4 Corollary 1). The most technical difficulty in studying such extensions is that we are going outside the category
of schemes. In fact an extension of a group scheme by an abelian scheme is
in particular a torsor under this abelian scheme and
in~\cite{R} XIII 3.2 Raynaud gives an example of a torsor under an abelian $S$-scheme which is not representable by a scheme. Nevertheless according to Lemma~\ref{algspace} such extensions of $S$-group schemes by abelian $S$-schemes are algebraic spaces over $S$.\\
We focus our attention on extensions of $S$-tori by abelian $S$-schemes.
Working over an arbitrary base scheme $S$ we prove that for such extensions it is equivalent to be a scheme or to be of finite order locally over $S$ for the Zariski topology. Then we show that
these extensions are in fact representable by schemes. Hence we can conclude that the extensions of $S$-tori by abelian $S$-schemes are of finite order locally over $S$.
We proceed in this way because it turns out to be easier to prove that these extensions are schemes
than to prove that they are of finite order locally over $S$, even though these two facts are equivalent in the end.\\
We start our investigation about extensions of tori by abelian schemes working over an algebraically closed field.

\begin{prop}\label{extTAcorps}
Let $k$ be an algebraically closed field. An extension $E$ of a
$k$-torus $T$ by an abelian $k$-variety $A$ is a connected smooth $k$-algebraic group. Moreover, if we denote by $A'$ and $T'$ respectively the abelian $k$-variety and the $k$-torus given by
Chevalley's decomposition (\cite{Ro} Theorem 16) of $E$, we have that
\begin{itemize}
    \item the torus $T'$ is the maximal torus of the extension $E$ and it is isogenous to the torus $T$;
    \item $\dim T' =\dim T$ and $\dim A' =\dim A$.
\end{itemize}
Moreover, $E$ is of finite order.
\end{prop}

\begin{proof} By Lemma~\ref{algspace} and by~\cite{A69} Lemma 4.2, the extension $E$ is
a connected smooth $k$-algebraic group.\\
Since the quotient $E /T'$ is an abelian variety, the torus $T'$ is the maximal torus of $E$ and so by~\cite{SGA3} Expos\'e XII Theorem 6.6 (d), via the surjective morphism $ E \rightarrow T$, the torus $T'$ goes onto $ T$. Moreover the kernel of this surjective morphism
 $T' \rightarrow T$ is $T'\cap A$ which is finite. Therefore the morphism
$T' \rightarrow T$ is an isogeny. This implies that $\dim T' =\dim T$ and that $\dim A' =\dim A$.\\
The existence of an isogeny between the two tori $T'$ and $T$ has a geometric implication: the extension $E$ is of finite order. In fact, let $n$ be a positive integer which annihilates the kernel of this isogeny.
By definition of push-down (see~\cite{S0} Chapter VII \S 1.1), the image $T''$ of the torus $T'$ in the push-down $n_*E$ of $E$ via the multiplication by $n$ on $A$, is isomorphic to the torus $T$ and this isomorphism between $T''$ and $T$ furnishes the section which splits the extension $n_*E.$
\end{proof}

Now we go back to the general case: let $S$ be an arbitrary scheme.

\begin{lem}\label{algspaceAT}
 An extension of an $S$-torus by an abelian $S$-scheme is an algebraic space over $S$ which is a group object and which is smooth (in particular flat), separated, of finite presentation over $S$ and with connected fibres.
\end{lem}

\begin{proof} Let $E$ be an extension of an $S$-torus $T$ by an abelian $S$-scheme $A$.
By Lemma~\ref{algspace}, the extension $E$ is an algebraic space over $S$ which is a group object. Clearly it has connected fibres.
Since the abelian scheme $A$ is proper and smooth over $S$, the $A$-torsor
$E$ is proper and smooth over $T$ and so $E$ is smooth, separated and of finite presentation over $S$.
\end{proof}

By Proposition~\ref{extTAcorps}, over each geometric point $\os$ of $S$ the  fibre $E_{\os}$ is a connected smooth algebraic group over the algebraically closed field $\oks$, and so we can generalized to the algebraic space $E$ the notion of abelian and reductive rank introduced by Grothendieck in~\cite{SGA3} Expos\'e X page 121:

\begin{itemize}
    \item the \textbf{abelian rank of $E$ at the point $s$}, denoted by $\rho_{\mathrm {ab}}(s),$ is the dimension of the
abelian $\oks$-variety appearing in Chevalley's decomposition of $E_{\os}$.
    \item the \textbf{reductive rank of $E$ at the point $s$}, denoted by $\rho_{\mathrm r}(s),$ is the dimension of the
maxinal tori of $E_{\os}$, i.e. the dimension of the $\oks$-torus appearing in Chevalley's decomposition of $E_{\os}$.
\end{itemize}

\begin{prop}\label{extTA-equiv}
Let $E$ be an extension of an $S$-torus by an abelian $S$-scheme.
The following conditions are equivalent:
\begin{description}
    \item[$(i)$] $E$ is a scheme,
    \item[$(ii)$] $E$ is of finite order locally over $S$ for the Zariski topology. If $S$ is quasi-compact, then the extension $E$ is globally of finite order.
\end{description}
\end{prop}

\begin{proof}  Let $E$ be an extension of an $S$-torus $T$ by an abelian $S$-scheme $A$.\\
$(i) \Rightarrow (ii) :$ Assume $E$ to be an $S$-scheme. Let $U$ be a quasi-compact open subset of $S$. We have to show that the extension $E_{U}$ of $T_{U}$ by $A_{U}$ is of finite order.
Denote by $\mathcal{T}$ (resp. by $\mathcal{MT}$) the functor of sub-tori  (resp. the functor of maximal sub-tori) of $E_{U}$ (see~\cite{SGA3} Expos\'e XV \S 8 for the definition of these functors). Since $E$ is a commutative group scheme the functor of sub-tori coincide with the functor of central sub-tori. In Lemma~\ref{algspaceAT} we have showed that
the extension $E_{U}$ is smooth, separated and of finite presentation over $U$. Moreover by Proposition~\ref{extTAcorps} its abelian rank and its reductive rank are locally constant functions over $U$. Therefore according to~\cite{SGA3} Expos\'e XV Corollary 8.11 and corollary 8.17, the functor $\mathcal{T}$ is representable by an \'etale and separated $U$-scheme and
the functor $\mathcal{MT}$ is representable by an open and closed sub-scheme of $\mathcal{T}$. In particular $\mathcal{MT}$ is \'etale over $\mathcal{T}$ and so over $U$. Since over each geometric point of $U$, the extension $E_{U}$ admits a unique maximal torus, by~\cite{EGAIV} 4 Corollary 17.9.5 $\mathcal{MT}$ is isomorphic to $U$, which implies that there is a unique maximal sub-torus $T'$ of $E_{U}$.\\
Over each geometric point of $U$, we have an epimorphism from the torus $T'$ to the torus $T_{U}$.
The kernel of this epimorphism is the scheme $T' \cap A_{U}$ which is flat and finite over $U$. Therefore $T' \rightarrow T_{U}$ is an isogeny.
Now let $n$ be a positive integer which annihilates the kernel of this isogeny. By definition of push-down, the image $T''$ of the torus $T'$ in the push-down $n_*E_{U}$ of $E$ via the multiplication by $n$ on $A_U$, is isomorphic to the torus $T_{U}$ and this isomorphism between $T''$ and $T_{U}$ furnishes the section which splits the extension $n_*E_{U}.$\\
$(ii) \Rightarrow (i):$ If the extension $E$ of $T$ by $A$ is trivial, then $E$ is isomorphic to the product $A \times T$ and so it is an $S$-scheme. If the extension $E$ of $T$ by $A$ is of order $n$, the extension $n_*E$ is trivial and hence, as we have just seen, it is an $S$-scheme. Consider the short exact sequence given by the multiplication by $n$
\[0   \longrightarrow E[n]  \longrightarrow  E  \longrightarrow n_*E \longrightarrow  0. \]
The kernel $E[n]$ of the multiplication by $n$ is an $A[n]$-torsor and so it is finite over $S$. By~\cite{A68} Theorem 3.3 $E[n]$ is then a scheme. Since the $S$-scheme $E[n]$ is affine over $S$, the $E[n]$-torsor $E$ is affine over the $S$-scheme $n_*E$ and therefore according to~\cite{SGA1} VIII Theorem 2.1  $E$ is an $S$-scheme.
\end{proof}

\begin{thm}\label{extTA}
Any extension of an $S$-torus by an abelian $S$-scheme is a scheme.
\end{thm}

\begin{proof}  Let $E$ be an extension of an $S$-torus $T$ by an abelian $S$-scheme $A$.
Since the question is local over $S$, we start doing two reduction steps:
\begin{itemize}
    \item by Lemma~\ref{algspaceAT} the algebraic space $E$ is of finite presentation over $S$ and so we can suppose $S$ to be an affine Noetherian scheme;
    \item by~\cite{A68} Theorem 3.2 an algebraic space $E$ is a scheme if and only if the reduced algebraic space $E_{\mathrm{red}}$ associated to $E$ is one. Hence we may assume $S$ to be reduced.
\end{itemize}
Let $s$ to be a maximal point of $S$, i.e. a generic point of an irreducible component of $S$. According to Proposition~\ref{extTAcorps}, the fibre $E_s$ over $s$ admits a maximal torus, which extends to a maximal torus of $E_U$, with $U$ an open non-empty subset of $S$ containing $s$. Doing the same thing at a maximal point of the complement of $U$ in $S$, by Noetherian recurrence on $S$, there exists a finite covering of $S$ by locally closed sub-schemes $S_i$,  for $i=1, \dots , r$, such that for each $i$ the restriction $E_i$ of the extension $E$ over $S_i$ admits a maximal torus $Z_i$. As we have observed in the proof of Proposition~\ref{extTA-equiv} $(i) \Rightarrow (ii)$, this torus $Z_i$ is isogenous to the restriction  $T_i$ of the torus $T$ over $S_i$ and if $n_i$ is a positive integer which annihilates the kernel of this isogeny, the push-down $(n_i)_*E_i$ of the extension $E_i$ via the multiplication by $n_i$ on the restriction of $A$ over $S_i$, is trivial since the image of the torus $Z_i$ in this push-down $(n_i)_*E_i$ is isomorphic to $T_i$.\\
Let $n$ be the least common multiple of $n_1, \dots, n_r$. If the extension $n_*E$ is a scheme, then $E$ is also one via the argument given in the proof of Proposition~\ref{extTA-equiv} $(ii) \Rightarrow (i).$
Therefore modulo multiplication by $n$, we can assume for each $i$ that the extension $E_i$ is trivial and the torus
  $Z_i$ is isomorphic to $T_i$. In particular, the tori $Z_i$ are locally closed subspaces of the extension $E$.\\
Denote by $Z$ the finite union of the tori $Z_i$ for $i=1, \dots , r$.
Since the tori $Z_i$ are locally closed in $E$ and $E$ is Noetherian, the set $Z$ is a globally constructible set of the underlying topological space of $E$. Set-theoretically $Z$ is characterized by the following property: if $x$ is a point of $E$ over the point $s$ of $S$, then $x$ is a point of $Z$ if and only if $x$ is a point of the maximal torus of the $\oks$-algebraic group $E_{\os}$.
This characterization of $Z$ shows that $Z$ commutes with base extensions, i.e. if $S'$ is a Noetherian $S$-scheme, if $E'$ is the algebraic space obtained from $E$ by base extension $S' \rightarrow S$, and if $Z'$ is the constructible set of $E'$ constructed as $Z$ in $E$, then $Z'$ is the inverse image of $Z$ by the canonical projection $pr:E'= E \times_{S} S' \rightarrow E$.\\
Now we will prove that this globally constructible set $Z$ is a closed subset of the algebraic space $E$. Using~\cite{LM-B} Corollary (5.9.3) we can generalize to algebraic spaces the characterization of closed constructible subsets of a scheme given in~\cite{EGAI} Chapter I Corollary (7.3.2): a constructible  subset of a Noetherian algebraic space is closed if and only if it is stable under specialization. In order to prove that $Z$ is closed under specialization we use the valuative criterion of specialization~\cite{LM-B} Proposition (7.2.1): if $x$ is a point of $Z$ and $y$ is a specialization of $x$, i.e. $y \in {\overline{\{x\}}}$, there is an $S$-scheme $L$, spectrum of a discrete valuation ring, and an $S$-morphism $\phi: L \rightarrow E$ which sends the generic point of $L$ to $x$ and the closed point of $L$ to $y$. Therefore, $Z$ is closed if and only if the image via $\phi$ of the closed point of $L$ is a point of $Z$. According to~\cite{An} Chap. IV Theorem 4.B the extension $E_L= E \times_{S} L$, which is obtained by base extension through the canonical morphism $L \rightarrow S$, is an $L$-scheme. Denote by $Z_L$ the inverse image of $Z$ by the canonical projection $E_L= E \times_{S} L \rightarrow E$. Using the set-theoretical characterization of $Z_L$, we observe that $Z_L$ is the maximal torus of the $L$-group scheme $E_L$ (see definition~\cite{SGA3} II Expos\'e XV 6.1) and hence it is closed.
The $S$-morphism $\phi: L \rightarrow E$ furnishes a section $\phi_L: L \rightarrow E_L$ of the $L$-scheme $E_L$, which sends the generic point of $L$ to $Z_L$. But since $Z_L$ is closed, also the closed point of $L$ maps to $Z_L$ via $\phi_L$. Therefore the morphism $\phi$, which is the composite of the section $\phi_L$ with the canonical projection $E_L \rightarrow E$, sends the closed point of $L$ to $Z$. This concludes the proof that $Z$ is a closed subset of $E$.\\
We endow $Z$ with the reduced induced closed algebraic subspace structure.
The algebraic subspace $Z$ is then proper over $T$.\\
If $s$ is  a point of $S$, according to Proposition~\ref{extTAcorps} the maximal torus of $E_{\overline s}$ is isogenous to $T_{\overline s}$. The set-theoretical characterization of $Z$ implies then that for each point $y$ of $T$ over $s$, the fibre $Z_y$ over $y$ is a scheme consisting of only one point which has the same residue field of the point $y$. Hence\begin{itemize}
    \item   by~\cite{LM-B}  Corollary (A.2.1), the morphism $Z \rightarrow T$ is finite. In particular, according to~\cite{A68} Theorem 3.3
the algebraic subspace $Z$ is a scheme.
    \item by~\cite{EGAI} Chapter I Proposition (3.7.1) (c) and Corollary (3.6.3) respectively, the morphism $Z \rightarrow T$ is a universal homeomorphism.
\end{itemize}
Denote by ${\overline E}$ the $A$-torsor over $Z$ obtained as pull-back of the $A$-torsor $E$ via the morphism $Z \rightarrow T$. The closed immersion $Z \rightarrow E$ trivialized this $A$-torsor ${\overline E}$, which is therefore a scheme. Because of the structure of $A$-torsor, the morphism $ {\overline E} \rightarrow E$ is a finite universal homeomorphism. In particular, if $U$ is an affine open subset of ${\overline E}$, there exists a Zariski open subset $V$ of $E$ such that $U$ is the inverse image of $V$ via this universal homeomorphism.
The restriction $U \rightarrow V$ of  $ {\overline E} \rightarrow E$ is again a finite universal homeomorphism and so, by Chevalley's Theorem~\cite{K} III 4.1 the Zariski open subset $V$ is affine. In this way we get an open affine covering of $E$, which is therefore a scheme.
\end{proof}

\begin{cor}\label{extTA-cor1}
Let $E$ be an extension of an $S$-torus $T$ by an abelian $S$-scheme $A$.
Then $E$ is of finite order locally over $S$ for the Zariski topology. If $S$ is quasi-compact, $E$ is globally of finite order. In particular, the $S$-sheaf ${\uExt}^1(T,A)$ is a torsion sheaf.
\end{cor}

\begin{cor}\label{extTA-cor3}
Let $A$ be an abelian $S$-scheme and let $T$ be
an $S$-torus.
The $S$-sheaf ${\uExt}^1(T,A)$ is an $S$-group scheme which is separated and \'etale over $S$.
\end{cor}

\begin{proof} Denote by ${\uExt}^1(T,A)[n]$ the kernel of the multiplication by $n$ on the sheaf ${\uExt}^1(T,A)$ for each integer $n$ bigger than 0. According to the above Corollary, the $S$-sheaf ${\uExt}^1(T,A)$ is a torsion sheaf and so it is the inductive limit of the
family of sheaves $\{{\uExt}^1(T,A)[n]\}_n$:
\begin{equation}\label{torsion-1}
    {\uExt}^1(T,A)= \lim_{\longrightarrow}  {\uExt}^1(T,A)[n].
\end{equation}
Consider the long exact sequence associated to the multiplication by $n$ on $T$:
\[
 \begin{array}{c}
0 \rightarrow {\uHom}(T,A) \rightarrow
{\uHom}(T,A) \rightarrow {\uHom}(T[n],A) \stackrel{d}{\rightarrow} {\uExt}^1(T,A)\\
\stackrel{n^*}{\rightarrow}  {\uExt}^1(T,A) \rightarrow {\uExt}^1(T[n],A)\\
\end{array}
\]
where $T[n]$ is the kernel of the multiplication by $n$ on $T$, $d$ is the connecting morphism, and $n^*$ is the pull-back of the extensions via the multiplication by $n$ on $T$. By Lemma~\ref{homTA}, the connecting morphism $d$ is injective. Therefore the $S$-sheaf ${\uHom}(T[n],A)$ is isomorphic to the $S$-sheaf ${\uExt}^1(T,A)[n]$ and the inductive limit~(\ref{torsion-1}) can be rewritten in the following way:
\begin{equation}\label{torsion-2}
    {\uExt}^1(T,A)= \lim_{\longrightarrow}  {\uHom}(T[n],A)
\end{equation}
where the inductive system is determined by the morphisms between the kernels
$T[n]$ given by the multiplication by integers: $m:T[mn] \rightarrow T[n].$
Since the multiplication by a nonzero integer on a torus is an epimorphism for the fppf topology, the morphisms ${\uHom}(T[n],A) \rightarrow {\uHom}(T[mn],A)$ of the inductive system are monomorphisms for the fppf topology. According to~\cite{EGAI} Chp. I (2.4.3) and Chp. 0 Proposition 4.5.4, in order to prove that the inductive limit~(\ref{torsion-2}) is a scheme, we have to verify that
\begin{description}
    \item[$(a)$] for each integer $n$, the $S$-sheaf ${\uHom}(T[n],A)$ is a scheme;
    \item[$(b)$] for each integer $m$, each monomorphism ${\uHom}(T[n],A) \rightarrow {\uHom}(T[mn],A)$ of the inductive system is an open immersion.
\end{description}
Condition $(a)$ is a consequence of~\cite{FGA} Expos\'e 221 4.c.\\
Before to prove $(b)$, we show that for each $n$ the scheme
${\uHom}(T[n],A[n])$ is \'etale over $S$.
In order to show this, it is enough to prove that the scheme ${\uHom}(A[n]^*,T[n]^*)$,
where $A[n]^*$ (resp. $T[n]^*$) is the Cartier dual of $A[n]$ (resp. $T[n]$), is \'etale over $S$. Consider the $S$-sheaf of morphisms ${\uMor}(A[n]^*,T[n]^*)$. Since the question is local over $S$, we can suppose $T[n]^*$ to be a constant group scheme and so for any $S$-scheme $S'$, to have a $S'$-morphism from $A[n]^*_{S'}$ to $T[n]^*_{S'}$ is equivalent to have a partition of $A[n]^*_{S'}$ in a finite number of open and closed subsets. By~\cite{EGAIV} 4 Lemma (18.5.3)
the $S$-sheaf ${\uMor}(A[n]^*,T[n]^*)$ is therefore an $S$-scheme which is affine, \'etale and of finite presentation over $S$. Consider now the morphism
\[ {\mathcal{H}}:{\uMor}(A[n]^*,T[n]^*) \longrightarrow {\uMor}(A[n]^* \times A[n]^*, T[n]^*) \]
which sends a morphism $f$ to the morphism $(x,y) \mapsto f(x) f(y) f(xy)^{-1}$. Since the scheme
 ${\uMor}(A[n]^* \times A[n]^*, T[n]^*)$ is \'etale, its unit section $\epsilon: S \rightarrow {\uMor}(A[n]^* \times A[n]^*, T[n]^*)$ is an open immersion.
Therefore the kernel $\mathcal{H}^{-1}(\epsilon(S))$ of $\mathcal{H}$, which is
${\uHom}(A[n]^*,T[n]^*)$, is an open subset of ${\uMor}(A[n]^*,T[n]^*)$ and so we can conclude that the scheme ${\uHom}(A[n]^*,T[n]^*)$ is \'etale over $S$.
Condition (b) is now a consequence of the fact that \'etale monomorphisms are open immersions.
\end{proof}

\subsection{The case in which only abelian schemes are involved}

We first study the $S$-sheaf of homomorphisms between two abelian schemes.

\begin{prop}\label{uhomAA}
Let $A$ and $B$ be two abelian $S$-schemes. Then the $S$-sheaf ${\uHom}(A,B)$ is an $S$-group scheme, which has constant geometric fibres defined by finitely generated free $\ZZ$-modules, and
 which is separated, unramified and essentially proper over $S$.
\end{prop}

\begin{proof} We can suppose the base scheme $S$ to be affine and Noetherian.
 By~\cite{A69} Corollary 6.2 (see also correctum~\cite{A74}
Appendix 1) the Hilbert functor ${\mathrm{Hilb}}_{(A \times_S B)/S}$ of
$A \times_S B$ is an algebraic space over $S$ which is locally of finite presentation and separated over $S$. The $S$-sheaf ${\uHom}(A,B)$ is then an algebraic subspace of
${\mathrm{Hilb}}_{(A \times_S B)/S}$ which is locally of finite presentation and separated over $S$. Moreover it has a group structure and it has constant geometric fibres defined by finitely generated free $\ZZ$-modules. Since its fibres are discrete, the structural morphism ${\uHom}(A,B) \rightarrow S$ is locally quasi-finite and so by~\cite{A68} Theorem 3.3 ${\uHom}(A,B)$ is a scheme. \\
In order to show that the scheme
${\uHom}(A,B)$ is unramified over $S,$ it is enough to show that if $S$ is the spectrum
of a local Artin ring $R$ with residue field $k$ and if $f: A \rightarrow B$ is an $S$-morphism whose restriction over ${\spec (k)}$ is trivial, then $f$ is trivial. Consider the restriction $f_{i}:A[q^i] \rightarrow B[q^i]$ of $f$ to the points of order $q^i$ with
$q$ an integer bigger than 1 which is coprime with the characteristic of $k$, and $i>0$. For each $i$, $f_{i}$ is a morphism of finite and \'etale $S$-schemes which is trivial over ${\spec (k)}$, and so it is trivial.
Since the family $(A[q^i])_i$ is schematically dense (see~\cite{SGA3} Expos\'e IX \S 4), we can conclude.\\
Finally, ${\uHom}(A,B)$ is essentially proper because of the Neronian property of abelian schemes over a discrete valuation ring.
\end{proof}

In order to investigate
extensions of abelian $S$-schemes by abelian $S$-schemes, we start working over a field:

\begin{prop}\label{extAAcorps}
Over a field,
\begin{enumerate}
    \item an extension of an abelian variety by an abelian variety is an abelian variety.
    \item an extension of an abelian variety by an abelian variety is of finite order.
\end{enumerate}
\end{prop}

\begin{proof}
 Let $E$ be an extension of an abelian $k$-variety $A$ by an abelian $k$-variety $B$.
By Lemma~\ref{algspace}, the extension $E$ is an algebraic space over $k$ which is a group object. Clearly $E$ is connected, proper and smooth over $k$ and so it is an abelian variety.\\
By Poincar\'e's complete reducibility theorem,
  there exists an abelian $k$-sub-variety $A'$ of $E$ such that $E=B \cdot A', $
$B \cap A' $ is finite and $\dim E= \dim B + \dim A'.$
We have a natural surjection
$A' \rightarrow A'/ B \cap A' =E/ B= A$. The kernel of this surjection is
$B \cap A' $ which is finite, and the surjection
$A'  \rightarrow A $ is an isogeny. This isogeny $A' \rightarrow A $ has a geometric consequence: the extension $E$ is of finite order. In fact
let $n$ be a positive integer $n$ which annihilates the kernel of this isogeny. By definition of push-down (see~\cite{S0} Chapter VII \S 1.1), the image $A''$ of the abelian variety $A'$ in the push-down $n_*E$ of $E$ via the multiplication by $n$ on $B$, is isomorphic to $A$ and this isomorphism between $A''$ and $A$ furnishes the section which split the extension $n_*E.$
\end{proof}

Now we go back to the general case: let $S$ be an arbitrary scheme.

\begin{lem}\label{algspaceAA}
An extension of an abelian $S$-scheme by an abelian $S$-scheme is an abelian $S$-scheme.
\end{lem}

\begin{proof} Let $E$ be an extension of an abelian $S$-scheme $A$ by an abelian $S$-scheme $B$. By Lemma~\ref{algspace}, the extension $E$ is an algebraic space over $S$ which is a group object. Clearly it has connected fibres. Moreover the $B$-torsor $E$ is proper and smooth
over $A$ and so $E$ is proper and smooth over $S$. Therefore by~\cite{F} Theorem 1.9 the extension $E$ is an abelian $S$-scheme.
\end{proof}

Therefore working with extensions of abelian $S$-schemes by abelian $S$-schemes we don't go outside the category of schemes. The difficulties with these extensions lie on the fact that it is not true that they are of finite order, not even of finite order locally over $S$ (see Remark below). Nevertheless, if we assume the base scheme $S$ to be integral and geometrically unibranched we will show that these extensions are of finite order locally over $S$. Recall that a scheme $S$ is \textbf{unibranched at the point} $\mathbf{s}$ if
the local ring ${\mathcal{O}}_{S,s}$ at the point $s$ is unibranched, i.e.
the ring $({\mathcal{O}}_{S,s})_{\mathrm{red}} = {\mathcal{O}}_{S,s} / {{\mathcal{N}}ilradical}$ is integral and the integral closure of $({\mathcal{O}}_{S,s})_{\mathrm{red}}$ is a local ring. A scheme $S$ is
\textbf{geometrically unibranched at the point} $\mathbf{s}$ if the local ring
${\mathcal{O}}_{S,s}$ at the point $s$ is unibranched and the residue field of the local ring, integral closure of $({\mathcal{O}}_{S,s})_{\mathrm{red}}$, is a radiciel extension of the residue field of ${\mathcal{O}}_{S,s}$.
Finally a scheme $S$ is \textbf{geometrically unibranched} if it is geometrically unibranched at each point $s$ of $S$. \\

\begin{rem}\label{extinfinite} \emph{(Extension of abelian schemes by abelian schemes which are of infinite order)} Let $X$ be the scheme obtained from infinitely many copies of the projective line ${\PP}_k^1$ over a field $k$  identifying the point 0 of a copy of ${\PP}_k^1$ with the point at the infinity of the following copy of ${\PP}_k^1$ and so on... The group $\ZZ$ acts on $X$ sending one copy of ${\PP}_k^1$ in the following one. The quotient $S=X / {\ZZ}$ is a scheme having a double point with distinct tangent lines. \\
Let $E$ be an elliptic curve over $S$. Consider the automorphism of $E \times_S E$ of infinite order
\begin{eqnarray}
\nonumber \mu: E \times_S E & \longrightarrow &  E \times_S E\\
\nonumber (a,b)  & \longmapsto & (a+b,b) .
\end{eqnarray}
 Consider also the trivial extension $E \times_S E $ of $E$ by $E$ over $S$.
 Since the automorphism $\mu$ respects exact sequences, we can use it in order to twist the cocycle defining the trivial extension, getting in this way an extension $\widetilde{E \times_S E}$ of $E$ by $E$ over $S$ which is of infinite order. Over $X$, the extension $\widetilde{E \times_S E}$ becomes the trivial extension.
\end{rem}

\begin{prop}\label{extAAtorsion}
Assume $S$ to be integral and geometrically unibranched. Let $E$ be an extension
of an abelian $S$-scheme by an abelian $S$-scheme.
Then $E$ is of finite order locally over $S$ for the Zariski topology. If $S$ is integral, geometrically unibranched and quasi-compact, $E$ is globally of finite order.
\end{prop}

\begin{proof}
Let $S= \spec(R)$ be an affine, integral and geometrically unibranched scheme.
 We can write $R$ as inductive limit of $\ZZ$-algebras $R_i$ of finite type. The $R_i$ are integral but not necessarily geometrically unibranched. Denote by $\eta_i$
 the generic point of $S_i= \spec(R_i)$.
  Let $E$ be an extension of an abelian $S$-scheme $A$ by an
abelian $S$-scheme $B$ and let $p:E \rightarrow A$ be the structural surjection of $E$ over $A$. Consider the descent $E_i$ of the extension $E$ over $S_i$ for $i$ big enough.
  According Proposition~\ref{extAAcorps} the extension $E_{i,\eta_i}$ is of finite order, say of order $n$. Hence there exists an open non-empty subset $U_i$ of $S_i$ and an $U_i$-morphism $q_i:A_{i,U_i} \rightarrow E_{i,U_i}$ such that the composite of $q_i$ with the restriction of $p_i$ over $U_i$ is the multiplication by $n$ (here $A_i$ and $p_i$ denote respectively the descent of $A$ and $p$ over $S_i$).
  Going back to $S$, there exists an open non-empty subset $U$ of $S$ and an $U$-morphism $q:A_{U} \rightarrow E_{U}$ such that the composite of $q$ with the restriction $p_U: E_U \rightarrow A_{U}$ of $p$ over $U$ is the multiplication by $n$: $p_U \circ q =n$. According to
Propositions~\ref{uhomAA}, the $S$-group scheme
${\uHom}(A,E)$ is unramified over $S$ and it satisfies the valuative criterion for properness. Since $S$ is integral and geometrically unibranched, by~\cite{EGAIV} 4 remark (18.10.20) the $U$-morphism $q$ extends to an  $S$-morphism $q:A \rightarrow E$ whose composition with the projection $p$ is the multiplication by $n$. By definition of pull-back, $q$ defines a section of the pull-back $n^*E$ of $E$ via the multiplication by $n$ on the abelian scheme $A$, which means that the extension $E$ is of order $n$.
\end{proof}

Let $A$ and $B$ abelian $S$-schemes. The $S$-sheaf ${\uExt}^1(A,B)$ is a complicated object: it is not representable by an $S$-scheme since it does not commute with
adic-limits (see condition ($F_3$) in~\cite{Mr}), and it is not a torsion sheaf as we have observed in remark~\ref{extinfinite}. Only if we assume the base scheme $S$ to be integral and geometrically unibranched,
using~\cite{EGAIV} 4 Theorem (18.10.1) we get that in the small \'etale site, the $S$-sheaf ${\uExt}^1(A,B)$ is a torsion sheaf.

\subsection{The case of extensions of abelian schemes by tori}

Through several ``d\'evissages'', using the results of the previous sections
we can study the group of homomorphisms, the $S$-sheaf of homomorphisms, the group of extensions and the $S$-sheaf of extensions involving
extensions of abelian schemes by tori.

\begin{prop}\label{ultimo}
Let $S$ be a scheme. Let $G_i$ (for $i=1,2$) be a commutative extension of an abelian
$S$-scheme $A_i$ by an $S$-torus $T_i$. The groups ${\Hom} (G_1,G_2)$ and ${\Ext}^1 (G_1,G_2)$ live respectively in the following diagrams
whose columns and rows are exact:
\[\begin{array}{ccccccccc}
& &   &  & 0 & & 0 & & {\Hom} (A_1,A_2) \\
& &   &  &  \downarrow & &\downarrow & & \downarrow \\
&  & 0 & \rightarrow & {\Hom} (G_1,T_2)  & \rightarrow & {\Hom} (T_1,T_2) &\rightarrow & {\Ext}^1 (A_1,T_2) \\
& & \downarrow  &  &  \downarrow & & \cong & & \downarrow  \\
0& \rightarrow& {\Hom} (A_1,G_2) & \rightarrow &  {\Hom} (G_1,G_2) & \rightarrow & {\Hom} (T_1,G_2)& \rightarrow & {\Ext}^1 (A_1,G_2)\\
\downarrow& & \downarrow  &  &  \downarrow & &\downarrow & & \downarrow \\
 0 & \rightarrow & {\Hom} (A_1,A_2) & \cong & {\Hom} (G_1,A_2)
& \rightarrow & 0 & \rightarrow &  {\Ext}^1 (A_1,A_2) \\
\downarrow & & \downarrow  &  &  \downarrow & &\downarrow & &  \\
 {\Hom} (T_1,T_2) & \rightarrow & {\Ext}^1 (A_1,T_2) & \rightarrow & {\Ext}^1 (G_1,T_2) & \rightarrow & {\Ext}^1 (T_1,T_2) & &   \\
\end{array}\]

\[\begin{array}{ccccccc}
0& \rightarrow & {\Hom} (T_1,T_2) & \cong & {\Hom} (T_1,G_2) & \rightarrow & 0  \\
\downarrow  & & \downarrow  &  &  \downarrow & &  \downarrow  \\
{\Hom} (A_1,A_2) & \rightarrow & {\Ext}^1 (A_1,T_2) & \rightarrow &  {\Ext}^1 (A_1,G_2) & \rightarrow & {\Ext}^1 (A_1,A_2)\\
 \cong & & \downarrow  &  &  \downarrow & & \downarrow \\
 {\Hom} (G_1,A_2) & \rightarrow & {\Ext}^1 (G_1,T_2) & \rightarrow & {\Ext}^1 (G_1,G_2) & \rightarrow &  {\Ext}^1 (G_1,A_2) \\
\downarrow & & \downarrow  &  &  \downarrow & & \downarrow  \\
0 & \rightarrow & {\Ext}^1 (T_1,T_2) & \rightarrow & {\Ext}^1 (T_1,G_2) & \rightarrow & {\Ext}^1 (T_1,A_2)    \\
\end{array}\]

\end{prop}

If $S$ is the spectrum of an algebraically closed field, we have that
\begin{itemize}
    \item the map ${\Ext}^1 (A_1,T_2) \rightarrow {\Ext}^1 (G_1,T_2)$ is surjective;
    \item the group ${\Ext}^1 (T_1,G_2)$ is a subgroup of ${\Ext}^1 (T_1,A_2).$ In particular the extensions of $T_1$ by $G_2$ are of finite order;
\item the extensions of $G_1$ by $A_2$ are of finite order.
\end{itemize}

For the $S$-sheaf of homomorphisms and the $S$-sheaf of extensions involving extensions of abelian schemes by tori the situation is quite similar: we get two diagrams which are analogous to the ones of Proposition~\ref{ultimo} and in particular as before we have that
\begin{itemize}
    \item the $S$-sheaf ${\uHom} (G_1,T_2)$ is a sub-sheaf of ${\uHom} (T_1,T_2);$
    \item the $S$-sheaf ${\uHom} (G_1,A_2)$ is isomorphic to the $S$-sheaf ${\uHom} (A_1,A_2);$
    \item the $S$-sheaf ${\uHom} (T_1,G_2)$ is isomorphic to the $S$-sheaf ${\uHom} (T_1,T_2);$
    \item the $S$-sheaf ${\uHom} (A_1,G_2)$ is a sub-sheaf of ${\uHom} (A_1,A_2).$
\end{itemize}
What is different with respect to the case of groups is that without putting any hypothesis on $S$ we have also that
\begin{itemize}
    \item the map ${\uExt}^1 (A_1,T_2) \rightarrow {\uExt}^1 (G_1,T_2)$ is surjective;
    \item the $S$-sheaf ${\uExt}^1 (T_1,G_2)$ is a sub-sheaf of ${\uExt}^1 (T_1,A_2),$ and so in particular it is a torsion sheaf.
\end{itemize}
If the base scheme $S$ is integral and geometrically unibrached, in the small \'etale site the $S$-sheaf ${\uExt}^1 (A_1,A_2) $
is a torsion sheaf and so
\begin{itemize}
    \item the $S$-sheaf ${\uExt}^1 (G_1,A_2)$ is a torsion sheaf in the small \'etale site.
\end{itemize}

%------------------------------------------------------------------------------
\section{Biextensions}

Let $S$ be a scheme.

Let $P,Q$ and $G$ be commutative $S$-group schemes. A \textbf{biextension of $(P,Q)$ by $G$} is a $G_{P\times Q}$-torsor $B$, endowed with a structure of commutative extension of $Q_P$
by $G_P$ and a structure of commutative extension of $P_Q$ by $G_Q,$ which are compatible one with another (for the definition of compatible extensions see~\cite{SGA7} Expos\'e VII Definition 2.1).
Let $B_i$ (for $i=1,2$) be a biextension of  $(P_i,Q_i)$ by $G_i,$ with
 $P_i,Q_i$ and $G_i$ commutative $S$-group schemes.
A \textbf{morphism of biextensions $B_1 \rightarrow B_2$} is a system $(F,f_{\scriptscriptstyle P},f_{\scriptscriptstyle Q},f_{\scriptscriptstyle G})$ where
$f_{\scriptscriptstyle P}:P_1 \rightarrow P_2,~ f_{\scriptscriptstyle Q}:Q_1 \rightarrow Q_2,~
 f_{\scriptscriptstyle G}:G_1 \rightarrow G_2$
are morphisms of $S$-group schemes, and $F:B_1 \rightarrow B_2$
is a morphism between the underlying sheaves of $B_1$ and $B_2$, such that $F$ is contemporaneously a morphism of extensions associated to the morphisms
\[f_{\scriptscriptstyle P}\times f_{\scriptscriptstyle Q}:Q_{1 \, P_1} \longrightarrow Q_{2 \,P_2},
\qquad f_{\scriptscriptstyle P}\times f_{\scriptscriptstyle G}:G_{1 \,P_1} \longrightarrow G_{2 \,P_2}\]
and a morphism of extensions associated to the morphisms
\[ f_{\scriptscriptstyle P}\times f_{\scriptscriptstyle Q}:P_{1 \,Q_1} \longrightarrow P_{2 \,Q_2},
\qquad f_{\scriptscriptstyle Q}\times f_{\scriptscriptstyle G}:G_{1 \,Q_1} \longrightarrow G_{2 \,Q_2}.\]
We denote by ${\bBiext}(P,Q;G)$ the category of biextensions of $(P,Q)$ by $G$. Let ${\Biext}^0(P,Q;G)$ be the group of automorphisms of any biextension
of $(P,Q)$ by $G$, and let ${\Biext}^1(P,Q;G)$ be the group of isomorphism
 classes of biextensions of $(P,Q)$ by $G$.
 By~\cite{SGA7} Expos\'e VII 2.7, the categories ${\bBiext}(P,Q;G)$ and ${\bBiext}(Q,P;G)$ are equivalent and so all what we prove for one of these categories is automatically proved also for the other.

Formalizing and generalizing Mumford's work~\cite{Mu1} on biextensions, in~\cite{SGA7} Expos\'e VII 3.6.5 and (3.7.4) Grothendieck points out the following homological interpretation of the groups ${\Biext}^0(P,Q;G)$ and ${\Biext}^1(P,Q;G)$:
\begin{equation}\label{homolointer}
 \begin{array}{c}
{\Biext}^0(P,Q;G) \cong {\Ext}^0 (P \stackrel{\scriptscriptstyle\LL}{\otimes} Q,G)\cong {\Hom}(P \otimes Q,G )\\
 {\Biext}^1(P,Q;G) \cong {\Ext}^1(P {\buildrel {\scriptscriptstyle \LL}
 \over \otimes}Q,G)\end{array}
 \end{equation}
where ${\Hom}(P \otimes Q,G)$ is the group of bilinear
 morphisms from $P\times Q$ to $G$,
 $P\stackrel{\scriptscriptstyle \LL}{\otimes}Q$ is the derived functor
 of the functor $Q \rightarrow P \otimes Q$ in the derived category of abelian sheaves. Using this homological interpretation of biextensions, in~\cite{SGA7} Expos\'e VIII (1.1.4) he gets the exact sequence of 5 terms
\begin{equation}\label{gen-5-terms}
 \begin{array}{c}
0 \rightarrow {\Ext}^1(P,{\uHom}(Q,G)) \rightarrow
{\Biext}^1(P,Q;G) \rightarrow {\Hom}(P,{\uExt}^1(Q,G))\\
\rightarrow {\Ext}^2(P,{\uHom}(Q,G))\rightarrow
{\Ext}^2(P, {\R}{\uHom}(Q,G)).\\
\end{array}
\end{equation}

We finish observing that since a biextension is in particular a torsor, by~\cite{LM-B} Corollary 10.4.2 we have

\begin{lem}
Any biextension of commutative $S$-group schemes is an algebraic space over $S$.
\end{lem}

\subsection{Biextensions involving locally constant group schemes}

\begin{prop}\label{biextGGX}
Let $X$ be an $S$-group scheme which is locally
for the \'etale topology a constant group scheme defined by a
 finitely generated free $\ZZ$-module, let $P$ be a smooth commutative $S$-group scheme with connected fibres, and let $Q$ be a commutative $S$-group scheme. Then,
\[{\bBiext}(P,Q;X)= {\bBiext}(Q,P;X)=0.\]
\end{prop}

\begin{proof} Since the $S$-sheaf ${\uHom}(P,X)$ is trivial,
\begin{itemize}
    \item the category $ {\bBiext}(Q,P;X)$ is rigid, and
    \item by the exact sequence of 5 terms (\ref{gen-5-terms}),
    there exists the canonical isomorphism
\[ {\Biext}^1(Q,P;X)\cong {\Hom}(Q, {\uExt}^1(P,X)).\]
\end{itemize}
Hence using again Proposition \ref{homextGX} we can conclude.
\end{proof}

Let $X$ be an $S$-group scheme which is locally
for the \'etale topology a constant group scheme defined by a
 finitely generated free $\ZZ$-module, and let $P$ be a commutative $S$-group scheme. The tensor product $P \otimes X$ is an algebraic space over $S$
 with the following property:
  there exist an \'etale surjective morphism $S' \rightarrow S$ for which
 $(P \otimes X)_{S'}$ is isomorphic to a finite product of copies of the $S'$-group scheme $P_{S'}.$ In particular, if $P$ is an $S$-group scheme which is locally
for the \'etale topology a constant group scheme defined by a
 finitely generated free $\ZZ$-module (resp. an $S$-torus, resp. an abelian $S$-scheme, resp. an extension of an abelian $S$-scheme by an $S$-torus), so is $P \otimes X$.

\begin{prop}\label{biextXGG}
Let $X$ be an $S$-group scheme which is locally
for the \'etale topology a constant group scheme defined by a
 finitely generated free $\ZZ$-module, and let $P_i$ (for $i=1,2$) be a commutative $S$-group scheme. Then,
\[{\bBiext}(X,P_1;P_2) \cong {\bExt}(X \otimes P_1,P_2)
\cong {\bBiext}(P_1,X;P_2) .\]
\end{prop}

\begin{proof} Consequence of~\cite{SGA7} Expos\'e VIII 1.2.
\end{proof}

\begin{cor}\label{biextXXG}
Let $X$ and $Y$ be $S$-group schemes which are locally
for the \'etale topology constant group schemes defined by
 finitely generated free $\ZZ$-modules and let $P$ be a commutative $S$-group scheme. Then,
\[{\bBiext}(X,Y;P)\cong {\bExt}(X \otimes Y,P).\]
Moreover, the objects of these categories are trivial locally over $S$ for the fppf to\-po\-lo\-gy.
\end{cor}

\begin{proof}
Consequence of the above Proposition and of Lemma~\ref{uextXP}.
\end{proof}

\subsection{Biextensions of group schemes by tori}

Biextensions of group schemes by tori have been studied by Grothendieck in~\cite{SGA7} Expos\'e VIII \S3. We recall here some of his results.

\begin{prop}\label{biextX?T}
Let $X$ be an $S$-group scheme which is locally
for the \'etale topology a constant group scheme defined by a
 finitely generated free $\ZZ$-module, let $A$ be an abelian $S$-scheme, and let $Y(1)$ and $Y'(1)$ be $S$-tori. Then,
\begin{enumerate}
    \item ${\bBiext}(X,Y'(1);Y(1)) \cong
{\bExt}(X \otimes Y'(1), Y(1))
\cong {\bBiext}(Y'(1),X;Y(1)). $\\
Moreover, the objects of these categories are trivial locally over $S$ for the fppf topology.
    \item ${\bBiext}(X,A;Y(1)) \cong {\bExt}(A, X^{\vee} \otimes Y(1))
\cong {\bBiext}(A,X;Y(1)).$\\
These categories are rigid and for the objects we have the canonical isomorphisms
\begin{eqnarray}
\nonumber{\Biext}^1(A,X;Y(1)) & \cong &{\mathrm{Ext}}^1(A,{\uHom}(X,Y(1)))\\
\nonumber {\Biext}^1(A,X;Y(1)) & \cong & {\Hom}(X,{\uExt}^1(A,Y(1))) ={\Hom}(X,Y \otimes A^*),
\end{eqnarray}
where $A^*$ is the dual abelian scheme ${\uExt}^1(A,{\ZZ}(1))$ of $A$.
    \item $ {\bBiext}(X,A;Y(1))\cong {\bExt}(X \otimes A , Y(1)) \cong  {\bBiext}(A,X;Y(1)).$\\
         In particular,
\[ {\bBiext}(X,A;Y(1)) \cong {\bExt}(X \otimes A, Y(1)) \cong {\bExt}(A, X^{\vee} \otimes Y(1)).\]
\end{enumerate}
\end{prop}

\begin{proof} (1) Consequence of Proposition~\ref{biextXGG}. The last assertion is due to~\cite{SGA7} Expos\'e VIII Proposition 3.3.1.\\
(2)~\cite{SGA7} Expos\'e VIII Proposition 3.7.\\
(3) Consequence of Proposition~\ref{biextXGG} and of assertion 2.
\end{proof}

\begin{prop}\label{biextGTT}
 Let $Y_1(1)$ and $Y_2(1)$ be $S$-tori and let $P$ be a commutative $S$-group scheme. Then,
\begin{equation}
 \label{biextGTT:1}{\bBiext}(P,Y_1(1);Y_2(1)) \cong {\bExt}(P, Y_1^{\vee} \otimes Y_2) \cong {\bBiext}(Y_1(1),P;Y_2(1)).
 \end{equation}
 where $Y_1^{\vee}$ is the character group of $Y_1(1)$ and $Y_2$ is the cocharacter group of $Y_2(1).$ Moreover,
 \begin{enumerate}
    \item If $P$ is a smooth commutative $S$-group scheme with connected fibres,
\[{\bBiext}(P,Y_1(1);Y_2(1)) = {\bBiext}(Y_1(1),P;Y_2(1))= 0.\]
    \item If $P$ is an $S$-group scheme $X$
which is locally for the \'etale topology a constant group scheme defined by a
finitely generated free $\ZZ$-module,
\[ {\bBiext}(X,Y_1(1);Y_2(1))\cong {\bExt}(X \otimes Y_1(1), Y_2(1))
\cong {\bExt}(X, Y_1^{\vee} \otimes Y_2).\]
 \end{enumerate}
\end{prop}

\begin{proof} By~\cite{SGA7} Expos\'e VIII Proposition 3.3.1 and paragraph 1.5, we have that
\[{\bBiext}(P,Y_1(1);Y_2(1)) \cong {\bExt}(P,{\uHom}(Y_1(1),Y_2(1))).\]
Using the canonical isomorphism  ${\uHom}(Y_1(1),Y_2(1)) \cong
 Y_1^{\vee} \otimes Y_2 ,$ we get the equivalence of categories~(\ref{biextGTT:1}).
 Assertion (1) is a consequence of Proposition~\ref{homextGX} and of~(\ref{biextGTT:1}). Assertion (2) is a consequence of Proposition~\ref{biextXGG} and of~(\ref{biextGTT:1}).
\end{proof}

The group of isomorphism classes of biextensions of abelian schemes by tori is isomorphic to the group of homomorphisms between abelian schemes. In fact,

\begin{prop}\label{biextAAT}
Let $A_i$ (for $i=1,2$) abelian $S$-schemes and let $Y(1)$ be a $S$-torus.
The category $ {\bBiext}(A_1,A_2;Y(1))$ of biextensions of $(A_1,A_2)$ by $Y(1)$ is rigid and for the objects we have the following canonical isomorphism
\[{\Biext}^1(A_1,A_2; Y(1)) \cong {\Hom}(A_1 , {\uExt}^1(A_2,Y(1))) = {\Hom}(A_1, Y \otimes A^*_2).\]
where $A^*_2$ is the dual abelian scheme ${\uExt}^1(A_2,{\ZZ}(1))$ of $A_2$.
\end{prop}

\begin{proof} Since the $S$-sheaf
${\uHom}(A,Y(1))$ is trivial,

\begin{itemize}
    \item the categories $ {\bBiext}(A_1,A_2;Y(1))$ is rigid, and
    \item by the exact sequence of 5 terms~(\ref{gen-5-terms}), there exists the canonical isomorphism
\[{\Biext}^1(A_1,A_2;Y(1))\cong {\Hom}(A_1 , {\uExt}^1(A_2,Y(1))).\]
\end{itemize}
\end{proof}

\subsection{Biextensions of group schemes by an abelian scheme}

Using the fact that the extensions of tori by abelian schemes are of finite order locally over $S$ (see Proposition~\ref{extTA-equiv} and Theorem~\ref{extTA}), we prove that there are no non-trivial biextensions of tori and divisible group schemes by abelian schemes.

\begin{thm}\label{biextPTA}
 Let $A$ be an abelian $S$-scheme, let $Y(1)$ be an $S$-torus and let $P$ be a divisible commutative $S$-group scheme locally of finite presentation over $S$, with connected fibres. Then,
\[{\bBiext}(P,Y(1);A)={\bBiext}(Y(1),P;A)=0\]
\end{thm}

\begin{proof} Since the $S$-sheaf ${\uHom}(Y(1),A)$ is trivial,
\begin{itemize}
    \item the categories ${\bBiext}(P,Y(1);A)$ is rigid, and
    \item by the exact sequence of 5 terms~(\ref{gen-5-terms}), there exists the canonical isomorphism
\[{\Biext}^1(P,Y(1);A) \cong {\Hom}(P,{\uExt}^1(Y(1),A)).\]
\end{itemize}
According to Corollary~\ref{extTA-cor3}, the $S$-sheaf ${\uExt}^1(Y(1),A)$ is an $S$-group scheme which is separated and \'etale over $S$. Since $P$ has connected fibres, by Lemma~\ref{conetal} we can conclude.
\end{proof}

There is a more down-to-earth proof of the above theorem: let $B$ be a biextension of $(P,Y(1))$ by $A$ and let
$\Phi:P \times Y(1) \times Y(1) \rightarrow A$ and $ \Psi:P \times P \times Y(1) \rightarrow A$ be the two co-cycles defining it.
 Since by Lemma~\ref{homTA} the group ${\Biext}^0(P,Y(1);A)$ is zero, in order to prove the triviality of the biextension $B$ we can work locally on the base and hence Corollary \ref{extTA-cor1} implies that it exists a positive integer $n$ such that $n \Phi=0.$
 The pull-back $(id \times n)^* B$ of $B$ via $id \times n:P \times Y(1) \rightarrow P \times Y(1)$ is defined by the trivial co-cycle $n\Phi:P \times Y(1) \times Y(1) \rightarrow A$ and by the co-cycle $ \Psi:P \times P \times Y(1) \rightarrow A$. But these 2 co-cycles are related one with another by the relation
\[
     \begin{array}{c}
     n\Phi (p_1+p_2,q_1,q_2)- n\Phi (p_1,q_1,q_2)-n \Phi (p_2,q_1,q_2)=\\
    \Psi (p_1,p_2,q_1+q_2)-\Psi (p_1,p_2,q_1)-\Psi (p_1,p_2,q_2)\\
    \end{array}
\]
for any $p_i$ point of $P$ and any $q_i$ point of $Y(1)$ (for $i=1,2$). Therefore by Lemma~\ref{homTA} also the co-cycle $\Psi$ is trivial.
 This means that
the pull-back $(id \times n)^* B$ of $B$ is trivial and so by the following Lemma, which is a consequence of~\cite{SGA7} Expos\'e VII Proposition 3.8.8, we can conclude.

\begin{lem}\label{biexttor}
Let $S$ be a scheme and let $P,Q$ and $G$ be three commutative $S$-group schemes.
Assume $P$ and $Q$ divisible. Then for integers $n$ and $m$, $n\not=0, m\not=0,$ the kernel of the map
\begin{eqnarray}
\nonumber \Upsilon_{(n \times m)}:
 {\Biext}^1(P,Q;G) & \longrightarrow & {\Biext}^1(P,Q;G) \\
\nonumber B & \longmapsto & (n \times m)^* B
\end{eqnarray}
sending each element of ${\Biext}^1(P,Q;G)$ to its pull-back via $n \times m:P \times Q \rightarrow P \times Q$ is trivial.
\end{lem}

We finish describing biextensions of locally constant group schemes and abelian schemes by abelian schemes and biextensions of locally constant schemes and tori by abelian schemes.

\begin{prop}\label{biextX?A}
Let $X$ be an $S$-group scheme which is locally
for the \'etale topology a constant group scheme defined by a
 finitely generated free $\ZZ$-module, let $Y(1)$ be an $S$-torus and let $A$ and $A'$ be abelian $S$-schemes. Then,
\begin{enumerate}
    \item ${\bBiext}(X,A';A) \cong {\bExt}(X \otimes A',A)
\cong {\bBiext}(A',X;A) .$\\
Moreover, if the base scheme $S$ is integral and geometrically unibranched, their objects are of finite order locally over $S$ for the Zariski topology.
    \item ${\bBiext}(X,Y(1);A) \cong {\bExt}(X \otimes Y(1),A)
\cong {\bBiext}(Y(1),X;A) .$\\
Moreover, these categories are rigid and their objects are of finite order locally over $S$ for the Zariski topology.
\end{enumerate}
\end{prop}

\begin{proof} (1) Consequence of Proposition~\ref{biextXGG} and of Proposition~\ref{extAAtorsion} \\
(2) The first assertion is a consequence of Proposition~\ref{biextXGG}. Lemma~\ref{homTA} implies that the categories $ {\bBiext}(X,Y(1);A)$ and ${\bExt}( X \otimes Y(1), A)$ are rigid. The last assertion is a consequence of Proposition~\ref{extTA-cor1}.
\end{proof}

\subsection{Biextensions of abelian schemes by abelian schemes}

Even though there are extensions of abelian schemes by abelian schemes which are of infinite order (see Remark~\ref{extinfinite}), in this section we show that there are no non-trivial biextensions of abelian $S$-schemes by abelian $S$-schemes. In order to do this we use the exact sequence of 5 terms~(\ref{gen-5-terms}) applied to tree abelian $S$-schemes $A,B$ and $C$
\begin{equation}\label{5-term}
    \begin{array}{c}
0 \rightarrow {\Ext}^1(A,{\uHom}(B,C)) \rightarrow
{\Biext}^1(A,B;C) \rightarrow {\Hom}(A,{\uExt}^1(B,C))\\
\rightarrow {\Ext}^2(A,{\uHom}(B,C))\rightarrow
{\Ext}^2(A, {\R}{\uHom}(B,C)),\\
\end{array}
\end{equation}
and we check that the first term ${\Ext}^1(A,{\uHom}(B,C))$ and the third term\\
${\Hom}(A,{\uExt}^1(B,C))$ are both trivial. The triviality of the first term is a consequence of Propositions~\ref{uhomAA} and~\ref{homextGX}.
Before we prove the triviality of the third term ${\Hom}(A,{\uExt}^1(B,C))$
 we need some Lemmas:

\begin{lem}\label{3-term-lem1}
Let $A, B$ and $C$ be three abelian $S$-schemes.
Then the group of morphisms from the abelian scheme $A$ to the $S$-sheaf
${\uExt}^1(B,C)$ is torsion-free.
\end{lem}

\begin{proof} Consequence of the fact that the abelian scheme $A$ is a divisible group scheme for the fppf topology.
\end{proof}

If $f:X \rightarrow S$ is a morphism of schemes and $\mathcal{G}$ is an abelian \'etale $X$-sheaf we have the Leray spectral sequence
\begin{equation}\label{leray:1}
 {\h}^p(S,{\R}^q f_* \mathcal{G})\Rightarrow  {\h}^{p+q}(X,\mathcal{G}).
\end{equation}
Recall also that a local ring $\mathcal{O}$, with residue field $k,$ is \textbf{henselian} if for each \'etale and separated morphism $S \rightarrow \spec(\mathcal{O})$, the set of $\spec(\mathcal{O})$-sections of $S$ is in one to one correspondence with the set of $\spec(k)$-sections of $S \times_{\spec(\mathcal{O})} \spec (k)$. A local ring is \textbf{strictly local} if it is an henselian ring and its residue field is separably closed.

\begin{lem}\label{3-term-lem1-5}
Let $S$ and $X$ be Noetherian schemes and let $f:X \rightarrow S$ be a morphism which is faithfully flat over $S$. Denote by $\eta$ the scheme of generic points of $X$ and $i: \eta \rightarrow X$ the canonical morphism. Let $\mathcal{H}$ be an abelian \'etale sheaf over $\eta$.Then
\begin{description}
    \item[$(a)$] the sheaves ${\R}^j i_*(\mathcal{H})$ and
     ${\R}^j (fi)_*(\mathcal{H})$ are torsion sheaves for each $j\geq 1$;
    \item[$(b)$] the sheaves ${\R}^j f_*(i_*(\mathcal{H}))$ are torsion sheaves for each $j\geq 1$.
\end{description}
\end{lem}

\begin{proof}
 $(a)$ We can assume $X$ (resp. $S$) to be the spectrum of a strictly local ring. We get assertion $(a)$ by applying the Leray spectral sequence~(\ref{leray:1})
 to the sheaf $\mathcal{H}$ and to the morphisms $i: \eta \rightarrow X$ and $fi: \eta \rightarrow S$ respectively, since $ {\h}^{j}(\eta,\mathcal{H})$ is a torsion group for each $j\geq 1$.\\
 $(b)$ Consequence of $(a)$.
\end{proof}

In the next proofs we need a very naive notion of constructible \'etale sheaf of $\ZZ$-modules and we are not really working in the framework of constructibility as in SGA4 but rather in a mild variant adapted to our needs. Assume $S$ to be a Noetherian scheme. Here an \'etale $S$-sheaf $\mathcal{F}$ of $\ZZ$-modules is said \textbf{constructible} if it exists a finite partition $S_i$ of $S$, with $S_i$ locally closed, such that $\mathcal{F}_{|S_i}$ is a locally constant sheaf whose fibres are finitely generated $\ZZ$-modules. This notion is stable by inverse image via morphisms of finite type and by direct image via finite morphisms.
In fact,\begin{itemize}
 \item If  $T \rightarrow S$ is a finite morphism of schemes, the function
$ s \mapsto \mathrm{separable~degree}\\ \mathrm{of~} T \times k(s)$ is constructible.
 \item If $S$ is a normal Noetherian irreducible scheme with fraction field $K$ and
if $T$ is the normalisation of $S$ in some finite \'etale Galois extension $L$ of $K$, with Galois group $G$, then  the class of conjugacy of the inertia group at a point $t$ of $T$,
varies in a constructible  way. (One needs such a fact to prove the constructibility of
the direct image of a sheaf concentrated on the generic point).
        \end{itemize}

\begin{lem}\label{3-term-lem2}
Let $S$ be a Noetherian scheme and let $f:X \rightarrow S$ be a morphism which is proper, smooth and with geometrically connected fibres. Moreover let $\mathcal{F}$ be an abelian \'etale constructible $S$-sheaf of $\ZZ$-modules. Denote by $\mathcal{K}$ the inverse image $f^* \mathcal{F}$ of $\mathcal{F}$ by the morphism $f$. Then
\begin{description}
    \item[$(a)$] the canonical morphism $\mathcal{F} \rightarrow f_*(
    \mathcal{K})$ is an isomorphism,
    \item[$(b)$] the $S$-sheaves ${\R}^j f_*(\mathcal{K})$ are torsion sheaves for each $j\geq 1$.
\end{description}
\end{lem}

\begin{proof} Remark that for the proof of this Lemma we cannot use the proper base change Theorem since $\mathcal{F}$ may not be a torsion sheaf.\\
 $(a)$ It suffices to compare the global sections over a strictly local base. This assertion is then a consequence of the fact that $f$ has a trivial Stein factorization.\\
$(b)$ We can suppose $S$ to be the spectrum of a strictly local Noetherian ring.
We prove $(b)$ by induction on the dimension of the support of $\mathcal{F}$.
Denote by $\zeta$ the scheme of generic points of $S$ and by $k: \zeta \rightarrow S$ the canonical morphism. Consider the canonical morphism
 $u: \mathcal{F} \rightarrow k_*(k^* \mathcal{F})$. The sheaves
 $k_*(k^* \mathcal{F}), \mathrm{Ker}(u)$ and $\mathrm{Coker}(u)$ are constructible sheaves over $S$. By induction, we can assume that  assertion $(b)$ is true for the sheaves $\mathrm{Ker}(u)$ and $\mathrm{Coker}(u)$. Therefore it is enough to prove this assertion for a sheaf
 $\mathcal{F}$ of the kind $k_*( \mathcal{G})$ where $\mathcal{G}$ is a sheaf over $\zeta$.\\
 Denote by $\eta$ the scheme of generic points of $X$ and by $i: \eta \rightarrow X$ the canonical morphism.
 Since $f$ is smooth, by ~\cite{SGA4} Expos\'e XVI Corollary 1.2 the sheaf $f^*(\mathcal{F})=f^*k_*(\mathcal{G})$ is the sheaf $i_*(\mathcal{H})$, where $\mathcal{H}$ the inverse image of $\mathcal{G}$ by the natural morphism $\eta \rightarrow \zeta.$
 Therefore we are reduced to prove the assertion $(b)$ for a sheaf over $X$ of the kind $i_*(\mathcal{H})$, but this is exactly Lemma~\ref{3-term-lem1-5}.
\end{proof}

\begin{lem}\label{3-term-lem3}
Let $S$ be a Noetherian scheme and let $G$ be an $S$-group scheme unramified over $S$. Denote by $G'$ the maximal \'etale open sub-scheme of $G$. Denote by $G_{\acute{e}t}$ (resp. $G_{fppf}$) the \'etale sheaf (resp. fppf sheaf) underlying $G$
on the small site over $S$. Idem for $G'$.
\begin{description}
    \item[$(a)$] We have that $G'_{\acute{e}t}=G_{\acute{e}t}, G'_{fppf}=G_{fppf}$ and the canonical applications
        ${\h}^i(S_{\acute{e}t},G_{\acute{e}t}) \rightarrow {\h}^i(S_{fppf},G_{fppf})$ are bijections;
    \item[$(b)$] Moreover, if $G$ is separated over $S$ and for each geometric point $\os$ of $S$ $G_{\os}$ is of finite type, then $G_{\acute{e}t}$ is a constructible sheaf of $\ZZ$-modules over $S$.
\end{description}
\end{lem}

\begin{proof} Remark that $G'$ is in fact a group sub-scheme of $G$ since it contains the unit section.\\
$(a)$ Consider a section of $G$ above $S$. We have to show that it factors through $G'$. Since $G$ is unramified over $S$, this section is an open immersion and so $G$ is flat over the points of this section.
Again by hypothesis of unramifiedness, $G$ is flat at a point $g$ if and only if $G$ is \'etale at $g$. This implies that $G$ is \'etale over the points of this section which factors therefore through $G'$.
Hence on the small site
over $S$, $G$ and $G'$ furnish the same \'etale sheaf (resp. fppf sheaf). Since $G'$ is smooth over $S$, the last assertion is a consequence of~\cite{G3} Theorem 11.7.\\
$(b)$ The hypothesis on $G$ are also satisfied by $G'$, and so we can assume $G=G'$. We are therefore reduce to show that for any integral sub-scheme $T$ of $S$, $G'_{T}$ is a locally constant scheme, whose fibres are finitely generated $\ZZ$-modules, over a non-empty open subset of $T$. Denote by $t$ the generic point of $T$. By hypothesis on the geometric fibres of $G$, there exists a finite \'etale extension $t' \rightarrow t$ which splits the fibre of $G$ over $t$. We can extend the morphism $t' \rightarrow t$ to a finite \'etale morphism $T' \rightarrow T$ such that $G_{T'}$ admits a constant group sub-scheme $H'$ which has the same fibre as $G$ over the generic point $t'$. The scheme $H'$ is then an open group sub-scheme of $G_{T'}$ which satisfies the valuative criterion for properness over $T'$ and therefore it is closed (since $G$ is separated over $S$.) Hence $H'= G_{T'}$.
\end{proof}

We will apply the two above Lemmas to a $S$-group scheme $G$ of the kind
${\uHom}(B,C)$, where $B$ and $C$ are two abelian $S$-schemes. Observe that in general the fibres of $G$ don't vary in a constructible way as the one of $G'$, i.e. $G_{\acute{e}t}$
is constructible, but not $G$ in general.

\begin{prop}\label{3-term}
Let $A, B$ and $C$ be three abelian $S$-schemes. Then,
\[ {\Hom}(A,{\uExt}^1(B,C)) =0 .\]
\end{prop}

\begin{proof} Since the question is local over $S$, by~\cite{EGAIV} 4 Proposition 18.8.18 and\\
\cite{EGAIV} 3 Theorem 8.8.2 we can suppose $S$ to be the spectrum of a strictly local Noetherian ring $\mathcal{O}$, with maximal ideal $\mathcal{M}$ and residue field $k$. For each $n\geq 1 $, denote by $S_n$ the $S$-scheme $\spec (\mathcal{O}/\mathcal{M}^n).$\\
In order to study the group ${\Hom}(A,{\uExt}^1(B,C))$
we used the exact sequence of 4 terms
\begin{equation}\label{4-term}
 \begin{array}{c}
0 \rightarrow {\h}^1(S,{\uHom}(B,C)) \rightarrow
{\Ext}^1(B,C) \rightarrow {\h}^0(S,{\uExt}^1(B,C))\\
\rightarrow {\h}^2(S,{\uHom}(B,C))
 \end{array}
\end{equation}
associated to the spectral sequence $ {\h}^p(S,{\uExt}^q(B,C))\Rightarrow {\Ext}^*(B,C),$ where the cohomology groups are computed with respect to the fppf topology.\\
Let $u$ be a morphism from $A$ to ${\uExt}^1(B,C)$.
We can consider the morphism $u: A \rightarrow {\uExt}^1(B,C)$
as a section of ${\uExt}^1(B,C)$ over $A$ satisfying some additive properties and the property that its restriction to the unit section of $A$ is trivial. The exact sequence~(\ref{4-term}) implies that the obstruction to represent $u$ as a global extension of $B_A$ by $C_A$ lies in the cohomology group ${\h}^2(A,{\uHom}(B_A,C_A))$. By Lemma~\ref{3-term-lem3} we can compute this cohomology group with respect to the \'etale topology
and we can reduce to a constructible sheaf. Since $S$ is the spectrum of a strictly local ring, the Leray spectral sequence~(\ref{leray:1}) implies that ${\h}^2(A,{\uHom}(B_A,C_A))$ is the group ${\h}^0(S, {\R}^2f_*{\uHom}(B_A,C_A))$, where $f:A \rightarrow S$ is the structural morphism of $A$. But then according to
Lemma~\ref{3-term-lem2} ${\h}^2(A,{\uHom}(B_A,C_A))$ is a torsion group annihilated by an integer, says $N$. Since by Lemma~\ref{3-term-lem1}
it is enough to prove that $N u=0,$ we can then suppose that $u$ comes from a global extension $E$ of $B_A$ by $C_A.$ The extension $E$ is an abelian $A$-scheme (see Lemma~\ref{algspaceAA}).
Since the restriction of $u$ over the unit section of $A$ is trivial, also the extension $E$ is trivializable over the unit section of $A$.\\
Now we will prove that for $n\geq 1$ the restriction $E_n$ of the extension E over the $S_n$-scheme $A_{S_n}$ is trivializable. In order to simplify the notations, we denote by $A_n$ the abelian $S_n$-scheme $A \times_S S_n$ and by $B_{A_n}$ (resp. $C_{A_n}$) the abelian $A_n$-scheme $B_A \times_A A_n$ (resp. $C_A \times_A A_n$).
\begin{itemize}
    \item We start with the case $n=1$.
Since $A_1$ is smooth over $S_1=\spec (k)$, $A_1$ is integral and geometrically unibranched. By Proposition~\ref{extAAtorsion} we have that the restriction $E_1$ of $E$ over $A_1$ is of finite order, which implies that the restriction of $u$ over $A_1$ is a torsion element and therefore by Lemma~\ref{3-term-lem1} it is trivial. Applying the exact sequence~(\ref{4-term}) to $S=A_1$, we observe that this restriction
of $u$ over $A_1$ comes from an element of ${\h}^1(A_1,{\uHom}(B_{A_1},C_{A_1})).$
But the group ${\h}^1(A_1,{\uHom}(B_{A_1},C_{A_1}))$ is trivial: in fact let $P$ be a
${\uHom}(B_{A_1},C_{A_1})$-torsor over $A_1$. If $\eta$ denote the generic point of $A_1$, the group  ${\h}^1(\eta,{\uHom}(B_{A_1},C_{A_1})_{\eta})$ is trivial (recall that $k$ is separably closed) and so $P$ has a section over the generic point $\eta$.
Since $A_1$ is integral and geometrically unibranched, by~\cite{EGAIV} 4 remark (18.10.20) this section over the generic point extends over the whole abelian variety $A_1$, i.e. $P$ is a trivial torsor.
Therefore we can conclude that the extension $E_1$ corresponding to the restriction of $u$ over $A_1$ is trivializable.
\item For each $n$ the restriction $E_n$ of $E$ over $A_n$ is trivializable over the unit section of $A_n$, that we identify with $S_n$. We have the following diagram
\[\begin{array}{ccccccccc}
  E_1 & \rightarrow & \cdots & \rightarrow  & E_n & \rightarrow  & E_{n+1}& \rightarrow &\cdots  \\
  \downarrow &  &  &  & \downarrow &  &\downarrow  &&\\
  A_1 & \rightarrow  & \cdots & \rightarrow  &A_n & \rightarrow  & A_{n+1} &\rightarrow  & \cdots \\
\end{array}\]
By~\cite{I} Proposition 3.2 (b) the liftings of the abelian scheme $E_n$ are classified by the group
${\Ext}^1(E_1, {\li}(E_1)^{\vee} \otimes \mathcal{M}^n/\mathcal{M}^{n+1}).$
According to the theory of universal vectorial extensions, we know that this group is isomorphic to the group
\begin{equation}\label{fine}
{\Hom}(\omega_{E_1}, {\li}(E_1)^{\vee} \otimes \mathcal{M}^n/
\mathcal{M}^{n+1})
\end{equation}
where $\omega_{E_1}$ is the vector group ${\Hom}({\Ext}^1(E_1,{\GG}_a),{\GG}_a)$. Therefore the liftings of $E_n$ are classified by morphisms between locally free $\mathcal{O}_{A_1}$-modules of finite rank. Since the $k$-abelian variety is projective, these morphisms are defines by constants. Moreover the restriction of $E_n$ over the unit section of $A_n$ is trivializable and so these constants are trivial, i.e. the group~(\ref{fine}) is trivial. This implies that there exists a unique way to lift the extension $E_n$ and therefore since the trivial way is one way to do it, the extension $E_n$ is trivializable.
\end{itemize}
 We know that the restriction of the extension $E$ over the unit section of $A$ and over $A_n$ for $n\geq 1$ is trivializable. We want to show that this global extension $E$ of $B_A$ by $C_A$ is zero as section of ${\uExt}^1(B,C)$ over $A$. Consider the subsheaf $\underline{\F}$ of ${\uHom}(B_A,E)$ consisting of those morphisms from the abelian
scheme $B_A$ to the abelian scheme $E$ which are the identity once composed with the projection $E \rightarrow B_A$: in other words, $\underline{\F}$ is the $A$-sheaf of the trivializations of the extension $E$. According to Proposition~\ref{uhomAA}, this sheaf $\underline{\F}$  is an $A$-scheme locally of finite presentation, separated and non-ramified over $A$.
The scheme ${\uHom}(B_A,C_A)$ acts freely and transitively on $\underline{\F}.$ We want to show that $\underline{\F}$ is in fact a
${\uHom}(B_A,C_A)$-torsor over the abelian scheme $A$. We have therefore to prove that the set of points of $\underline{\F}$ where the structural morphism $\underline{\F} \rightarrow A$ is \'etale is sent surjectively to $A$. According to~\cite{EGAIV} 4 Theorem (17.6.1) it is enough to check that the set of points of $\underline{\F}$ where $\underline{\F} \rightarrow A$ is flat is sent surjectively to $A$.
Choose a trivialization $\sigma$ of the extension $E$ over the unit section of the abelian scheme $A$. Since the extension $E_n$ over $A_n$ is trivializable, there exists a section of the scheme $\underline{\F}$ over $A_n$. Modulo translation, we can suppose that these section coincides with the restriction of $\sigma$ over the unit section of the abelian scheme $A_n$. In this way we get a well-defined section
$\tau_n$ the scheme $\underline{\F}$ over $A_n$.
 In other words, we have used the trivialization $\sigma$ in order to find a ``compatible'' family of trivialization $\{\tau_n\}_n$. In particular $\underline{\F} \times_A A_n \rightarrow A_n$ is \'etale along $\tau_n$. Taking the limit over $n$, we get that the structural morphism $\underline{\F} \rightarrow A$ is \'etale along $\tau_n$ for each $n$. Hence if we denote by $U$ the open subset of $\underline{\F}$ where $\underline{\F}$ is \'etale over $A$, we have that $U$ contains $\tau_n$ for each $n$. Since the restriction of the structural morphism $\underline{\F} \rightarrow A$ to $U$ is locally of finite presentation and flat, by~\cite{EGAIV} 2 Theorem (2.4.6) it is universally open and so the image of $U$ in $A$ is an open subset $V$ of $A$. The open subset $U$ contains $\tau_1$ and therefore the open subset $V$ of $A$ contains the closed fibre $A_1$. But the abelian scheme $A$ is proper (hence universally closed) over $S$ and so $V$ is equal to $A$. This finishes the proof that $\underline{\F} $ is a ${\uHom}(B_A,C_A)$-torsor.  \\
The torsor $\underline{\F}$ is an element of ${\h}^1(A,{\uHom}(B_A,C_A))$ whose image in ${\Ext}^1(B_A,C_A)$ via the exact sequence~(\ref{4-term}) is the global extension $E$. Therefore the image of this extension $E$
in ${\h}^0(A,{\uExt}^1(B_A,C_A))$ via the exact sequence~(\ref{4-term}) is zero, i.e. this extension $E$ of $B_A$ by $C_A$ is zero as section of ${\uExt}^1(B,C)$ over $A$.
But this section was represented by the morphism $u: A \rightarrow {\uExt}^1(B,C)$ and so $u$ is trivial.
\end{proof}

\begin{thm}\label{biextAAA}
Let $A_i$ (for $i=1,2,3$) be an abelian $S$-scheme. Then,
\[{\bBiext}(A_1,A_2;A_3)=0.\]
\end{thm}

\begin{proof} Since
\[{\Biext}^0(A_1,A_2;A_3) \cong {\Hom}(A_1 \otimes A_2,A_3) \cong {\Hom}(A_1, {\uHom}(A_2,A_3)),\]
Proposition~\ref{uhomAA} and Lemma~\ref{conetal} imply that
the category ${\bBiext}(A_1,A_2;A_3)$ is rigid. Using Propositions~\ref{uhomAA},~\ref{homextGX} and~\ref{3-term}, from the exact sequence of 5 terms~(\ref{5-term})
 we get that the group ${\Biext}^1(A_1,A_2;A_3)$
is trivial.
\end{proof}

\subsection{Biextensions of extensions of abelian schemes by tori}

Through several ``d\'evissages'', using Theorems~\ref{biextPTA} and~\ref{biextAAA} we prove now the main theorem of this paper. We start with a first ``d\'evissage'':

\begin{prop}\label{biextGGA}
 Let $A$ be an abelian $S$-scheme and let $G_i$ (for $i=1,2$)
be a commutative extension of an abelian $S$-scheme $A_i$ by an $S$-torus $Y_i(1)$. Then,
\[ {\bBiext}(G_1,G_2;A) =0\]
\end{prop}

\begin{proof}
According to the homological interpretation~(\ref{homolointer})
 of the groups ${\Biext}^i$ (for $i=0,1$), from the short exact sequence
 $0 \rightarrow Y_2(1) \rightarrow G_2 \rightarrow A_2 \rightarrow 0$,
 we have the two long exact sequences
 \[ \begin{array}{c}
 0 \rightarrow {\Biext}^0 (A_1,A_2;A)\rightarrow {\Biext}^0 (A_1,G_2;A)\rightarrow {\Biext}^0 (A_1,Y_2(1);A) \rightarrow \\
\rightarrow {\Biext}^1 (A_1,A_2;A)\rightarrow {\Biext}^1 (A_1,G_2;A)\rightarrow
{\Biext}^1(A_1,Y_2(1);A) \rightarrow ...\\
\end{array}\]
\[ \begin{array}{c}
0 \rightarrow {\Biext}^0 (Y_1(1),A_2;A)\rightarrow {\Biext}^0 (Y_1(1),G_2;A)\rightarrow {\Biext}^0 (Y_1(1),Y_2(1);A) \rightarrow \\
\rightarrow {\Biext}^1 (Y_1(1),A_2;A)\rightarrow {\Biext}^1 (Y_1(1),G_2;A)\rightarrow
{\Biext}^1(Y_1(1),Y_2(1);A) \rightarrow ...\\
\end{array}\]
By Theorem \ref{biextPTA} and Theorem \ref{biextAAA}, these long exact sequences furnish the relations
\[{\bBiext}(A_1,G_2;A)=0  \qquad  \qquad {\bBiext}(Y_1(1),G_2;A)=0. \]
On the other hand, from the exact sequence
 $0 \rightarrow Y_1(1) \rightarrow G_1 \rightarrow A_1 \rightarrow 0$,
 we get the long exact sequence
\[\begin{array}{c}
 0 \rightarrow {\Biext}^0 (A_1,G_2;A)\rightarrow {\Biext}^0 (G_1,G_2;A)\rightarrow
{\Biext}^0 (Y_1(1),G_2;A) \rightarrow \\
\rightarrow {\Biext}^1 (A_1,G_2;A)\rightarrow {\Biext}^1 (G_1,G_2;A)\rightarrow
{\Biext}^1(Y_1(1),G_2;A) \rightarrow .....\\
\end{array}\]
Using this long exact sequence and the above relations, we can conclude.
\end{proof}

\begin{thm}\label{mthm1}
Let $S$ be a scheme.
Let $G_i$ (for $i=1,2,3$) be a commutative extension of an abelian
$S$-scheme $A_i$ by an $S$-torus $Y_i(1)$. In the topos $\mathrm{{\mathbf{T}_{fppf}}}$,
the category of biextensions of $(G_1,G_2)$ by $G_3$ is equivalent
to the category of biextensions of the underlying abelian $S$-schemes $(A_1,A_2)$ by the underlying $S$-torus $Y_3(1)$:
$$ {\bBiext}(G_1,G_2;G_3) \cong {\bBiext}(A_1,A_2;Y_3(1)) $$
\pn In particular, for $i=0,1$, we have the isomorphisms
\[{\Biext}^i(G_1,G_2;G_3)\cong {\Biext}^i(A_1,A_2;Y_3(1)).\]
\end{thm}

\begin{proof}
We will prove the following equivalences of categories:
\begin{equation}\label{biextGGG}
\begin{array}{c}
 {\bBiext}(G_1,G_2;Y_3(1)) \cong {\bBiext}(A_1,A_2;Y_3(1)) \\
 {\bBiext}(G_1,G_2;G_3) \cong {\bBiext}(G_1,G_2;Y_3(1))\\
\end{array}
\end{equation}
By~\cite{SGA3} Expos\'e X Corollary 4.5, we can suppose
that the tori are split (if necessary we localize over $S$ for the \'etale topology) and therefore we can assume that $Y_3(1)$ is ${\GG}_m^{\rk Y_3}.$
Since the categories ${\bBiext}(G_1,G_2;{\GG}_m)$ and ${\bBiext}(A_1,A_2;{\GG}_m)$
are additive in the variable ${\GG}_m$ (cf.~\cite{SGA7} I Expos\'e VII (2.4.2)), for the first equivalence of categories of (\ref{biextGGG}) it is enough to prove that
 \[{\bBiext}(G_1,G_2;{\GG}_m)\cong {\bBiext}(A_1,A_2;{\GG}_m)\]
and this is done in~\cite{SGA7} Expos\'e VIII (3.6.1). \\
By the homological interpretation~(\ref{homolointer}) of the groups ${\Biext}^i$ (for $i=0,1$), from the short exact sequence
 $0 \rightarrow Y_3(1) \rightarrow G_3 \rightarrow A_3 \rightarrow 0$
 we have the long exact sequence
\[\begin{array}{c}
 0 \rightarrow {\Biext}^0 (G_1,G_2;Y_3(1))\rightarrow {\Biext}^0 (G_1,G_2;G_3)
\rightarrow {\Biext}^0 (G_1,G_2;A_3) \rightarrow \\
\rightarrow {\Biext}^1 (G_1,G_2;Y_3(1))\rightarrow {\Biext}^1 (G_1,G_2;G_3)
\rightarrow {\Biext}^1(G_1,G_2;A_3) \rightarrow .....\\
\end{array}\]
Using Proposition~\ref{biextGGA}, we get the second equivalence of categories of (\ref{biextGGG}).
\end{proof}

\begin{rem} The above Theorem says essentially that each biextension
$\mathcal{B}$ of $(G_1, G_2)$ by $G_3$ comes from a biextension $B$
 of the underlying abelian $S$-schemes $(A_1,A_2)$ by the underlying $S$-torus $Y_3(1)$.
 If for $i=1,2,3,$ we denote by $\pi_i:G_i \rightarrow A_i$
 the projection of $G_i$ over $A_i$ and by $\iota_i:Y_i(1) \rightarrow
G_i$ the inclusion of $Y_i(1)$ in $G_i,$ we can describe explicitly the biextension
$\mathcal{B}$ of $(G_1, G_2)$ by $G_3$ in term of the corresponding biextension $B$
  of $(A_1, A_2)$ by $Y_3(1)$ as follow: $\mathcal{B}$ is the push-down
by $\iota_{3}$ of the biextension of $(G_1, G_2)$ by $Y_3(1)$
which is the pull-back by $(\pi_1,\pi_2)$ of $B$, i.e. $\mathcal{B} = \iota_{3\,*}(\pi_1,\pi_2)^*B.$
\end{rem}

% ------------------------------------------------------------------------

% ---------------------------------------------------------------------------

\end{document}